\documentclass[12pt]{article}
\usepackage{amsmath,amssymb,amsthm}
\usepackage[margin=1in]{geometry}
\usepackage{booktabs}

\newtheorem{thm}{Theorem}[section]
\newtheorem{lem}[thm]{Lemma}
\newtheorem{prop}[thm]{Proposition}

\theoremstyle{definition}
\newtheorem{defn}[thm]{Definition}
\theoremstyle{remark}
\newtheorem{rem}[thm]{Remark}

\newcommand{\Z}{\mathbb{Z}}
\newcommand{\E}{\mathbb{E}}
\newcommand{\Prob}{\mathbb{P}}

\title{A Faber--Krahn theorem for distance graphs on the integers}
\author{Suil O\thanks{Department of Applied Mathematics and Statistics, The State University of New York, Korea, Incheon, 21985, suil.o@sunykorea.ac.kr. Research supported by the National Research Foundation of Korea (NRF) grant funded by the Korea government(MSIT) No. RS-2025-23523950.}}
\date{}

\begin{document}
\maketitle

\begin{abstract}
Let $D$ be a finite set of positive integers with $|D|\ge2$ and $\gcd(D)=1$, let $\omega_k>0$ for $k\in D$, and consider the random walk on $\Z$ whose steps $\pm k$ ($k\in D$) have probabilities proportional to $\omega_k$. For a set $S\subseteq\Z$ that is bounded below and contains a half-line, the nonnegative function that vanishes off $S$, is harmonic on $S$ and differs from the identity by a bounded function satisfies $h_S(i)=i+\beta(S)+o(1)$ as $i\to\infty$; we call $\beta(S)$ the offset of $S$. In this paper, we prove that among the sets that contain $[N,\infty)$ and have a prescribed number of elements below $N$, the half-line has the largest offset, uniquely and with an explicit gap. We apply this to the weighted distance graph $G(\Z,\omega)$, in which $i$ and $j$ are joined by an edge of weight $\omega_{|i-j|}$ whenever $|i-j|\in D$. For every sufficiently large $m$, the block $\{0,1,\dots,m-1\}$ is, up to translation, the unique $m$-element set whose induced subgraph has the largest spectral radius, or equivalently the smallest first Dirichlet eigenvalue; this is a Faber--Krahn theorem for $G(\Z,\omega)$. For the unweighted distance graphs with $D=\{2,3\}$, $\{2,3,4\}$ and $\{3,4,5\}$, a computer-assisted proof shows that the block is extremal for every $m\ge7$, $8$ and $10$, respectively, uniquely except for $D=\{3,4,5\}$ and $m=10$, where one other set, up to reflection, induces an isomorphic graph.


\end{abstract}

\section{Introduction}

For a finite set $D\subseteq\Z_{>0}$, the \emph{distance graph} $G(\Z,D)$ has vertex set $\Z$ and $i\sim j$ if and only if $|i-j|\in D$. These graphs were introduced by Eggleton, Erd\H{o}s and Skilton~\cite{EES} and are usually studied through their coloring parameters~\cite{LZ}. Here we study their finite induced subgraphs through the spectral radius, and we allow weights. Let $\omega\colon D\to(0,\infty)$, and put $\omega_{-k}=\omega_k$ for $k\in D$ and $\omega_k=0$ for $|k|\notin D$. The \emph{weighted distance graph} $G(\Z,\omega)$ has adjacency matrix $A$ with $A_{ij}=\omega_{i-j}$, so every vertex has weighted degree $W=2\sum_{k\in D}\omega_k$. For finite $S\subseteq\Z$ let $\rho(S)$ be the spectral radius of the principal submatrix $A_S$, and put
\[
\rho_m(\omega)=\max\{\rho(S): S\subseteq\Z,\ |S|=m\}.
\]
The distance graph $G(\Z,D)$ is the case $\omega\equiv1$, and we then write $\rho_m(D)$. A set $S$ with $|S|=m$ and $\rho(S)=\rho_m(\omega)$ is \emph{extremal}. We write $B_m=\{0,1,\dots,m-1\}$ for the \emph{block}, and $[a,b]=\{a,a+1,\dots,b\}$ for an integer interval. For graph-theoretic and spectral terminology we follow~\cite{West}, \cite{BH} and~\cite{GR}.

The quantity $W-\rho(S)$ is the first eigenvalue of the Dirichlet Laplacian $WI-A_S$ of $S$ in $G(\Z,\omega)$, so $\rho_m(\omega)$ asks for the $m$-sets with the smallest first Dirichlet eigenvalue: it is the \emph{Faber--Krahn problem} on $G(\Z,\omega)$. On graphs this problem was raised by Friedman~\cite{Friedman}. For regular trees it was solved by Leydold~\cite{Leydold} and Pruss~\cite{Pruss}, for trees with a given degree sequence by B\i y\i ko\u{g}lu and Leydold~\cite{BL}, and for trees with a given matching number by Lin, Liu and You~\cite{LLY}. On the nearest-neighbor lattice $\Z^d$ with $d\ge2$ no exact description of the minimizers is known; they are close to discrete balls for large $m$~\cite{SR,CKLM}. On $\Z$ with $D=\{1\}$ the problem is trivial, since the induced subgraphs are disjoint unions of paths and the path $B_m$ wins. The graph $G(\Z,\omega)$ is the one-dimensional lattice with the steps in $D$ allowed, and the question is whether the interval still wins once the steps are not nearest-neighbor. It does not always win. If $D=\{d\}$, then the extremal sets are arithmetic progressions with difference $d$, and if $\gcd(D)=d\ge2$, then $G(\Z,D)$ is a disjoint union of copies of $G(\Z,D/d)$ and the extremal sets are not intervals. Even when $\gcd(D)=1$ the interval loses for small $m$. By exhaustive search, for $D=\{2,3,4\}$ and $3\le m\le6$ the unique extremal set is $\{0,2,\dots,2(m-1)\}$, which induces the square of a path, and for $D=\{3,4,5\}$ and $m=6$ it is $\{0,1,4,5,8,9\}$, a union of pairs placed with period $4$. As $m$ grows the runs lengthen; for $\{2,3,4\}$ and $m=7$ the unique extremal set is $\{0\}\cup[2,6]\cup\{8\}$, and for $\{3,4,5\}$ and $m=8$ it is $\{0,3,4,5,7,8,9,12\}$. In this paper, we prove that the block takes over for every $D$ with $|D|\ge2$ and $\gcd(D)=1$ and for every weight on $D$, and we determine the exact threshold for three distance graphs.

\begin{thm}\label{thm:all}
\begin{enumerate}
\item Let $D$ be a finite set of positive integers with $|D|\ge2$ and $\gcd(D)=1$, and let $\omega\colon D\to(0,\infty)$. There is $m_0(\omega)$ such that for every $m\ge m_0(\omega)$ the block $B_m$ is, up to translation, the unique extremal set for $\rho_m(\omega)$.
\item Let $D\in\{\{2,3\},\{2,3,4\},\{3,4,5\}\}$ and $m_0=7$, $8$, $10$ respectively. For every $m\ge m_0$ we have $\rho_m(D)=\rho(B_m)$. For $m\ge m_0$, and $m\ge11$ when $D=\{3,4,5\}$, the block $B_m$ is the unique extremal set up to translation. For $D=\{3,4,5\}$ and $m=10$ the extremal sets are $B_{10}$ and $[0,10]\setminus\{2\}$ up to translation and reflection, and they induce isomorphic graphs.
\end{enumerate}
\end{thm}

The thresholds in part~(2) are sharp: for $m=m_0-1$ the sets $\{0,2,3,4,5,7\}$, $\{0\}\cup[2,6]\cup\{8\}$ and $\{0,1,3,4,5,7,8,9,12\}$, respectively, have larger spectral radius than the block.

We first explain why the block wins. Since $W\|x\|^2-x^\top Ax=\sum_{k\in D}\omega_k\sum_i(x_i-x_{i+k})^2$, the number $W-\rho(S)$ is the least Dirichlet energy of a unit vector supported on $S$, and since $\gcd(D)=1$ the energy is small only for vectors that vary slowly. Such a vector sees $S$ through its \emph{effective length}: the number of vertices, corrected by an offset at each end. A competitor of the block can lose effective length in two ways. A hole away from the ends wastes a position where the Perron vector is large: the potential at such a hole is large, and moving the vertex with the smallest Perron entry into it increases the spectral radius (Proposition~\ref{prop:hole}). A flux argument turns this into the statement that every hole of an extremal set lies within bounded distance of an end (Proposition~\ref{prop:bulk}). A hole near an end changes the offset of that end, and among all ends with the same number of vertices the contiguous end has the largest offset. This last statement is the heart of the proof, and it is a statement about random walks.

Let $(X_t)$ be the random walk on $\Z$ with transition probabilities $p(x,y)=\omega_{x-y}/W$, so that the steps $\pm k$ have probabilities proportional to $\omega_k$. Call $S\subseteq\Z$ \emph{admissible} if it is bounded below and contains a half-line $[N,\infty)$. For such $S$ there is a unique nonnegative function $h_S$ that vanishes off $S$, is harmonic on $S$ for the walk, and satisfies $h_S(i)=i+\beta(S)+o(1)$ as $i\to\infty$ (Lemma~\ref{lem:h}); the number $\beta(S)$ is the \emph{offset} of $S$. Seen from far away, $S$ acts like the shifted half-line $[-\beta(S),\infty)$. Let $\sigma^2=\sum_{k\in D}\omega_kk^2\big/\sum_{k\in D}\omega_k$ be the variance of a step, $h_0=h_{[0,\infty)}$, and $G_0(x,y)$ the expected number of visits to $y$ of the walk started at $x$ before it leaves $[0,\infty)$.

\begin{thm}\label{thm:main}
Let $D$ be a finite set of positive integers with $|D|\ge2$ and $\gcd(D)=1$, and let $\omega\colon D\to(0,\infty)$. Let $S$ be admissible with $S\supseteq[N,\infty)$, and let $c=|S\cap(-\infty,N)|$. Then
\[
\beta(S)\le\beta([N-c,\infty)),
\]
with equality if and only if $S=[N-c,\infty)$. Moreover, if $S$ is not a half-line, then
\[
\beta([N-c,\infty))-\beta(S)\ \ge\ \gamma(\omega):=\frac{2}{\sigma^2}\cdot\frac{h_0(1)^2}{G_0(1,1)}-1,
\]
with equality if and only if $S$ is a translate of $\{-1\}\cup[1,\infty)$, and $\gamma(\omega)>0$.
\end{thm}

At positive frequency the offset becomes a phase in the quantization condition of an approximate eigenvector, so a non-contiguous end shortens the effective length by at least $\gamma(\omega)$ and lowers the spectral radius by order $m^{-3}$ (Proposition~\ref{prop:ends}). The block induces a banded Toeplitz matrix, whose extreme eigenvalue is well understood~\cite{BG}; the comparison needs the next term, of order $m^{-3}$, and that is where the offsets enter.

The proof of Theorem~\ref{thm:all}(1) gives no explicit $m_0(\omega)$; see Section~\ref{sec:remarks}. For part~(2) the same two mechanisms are made effective by an induction on $m$: the induction hypothesis bounds the edge boundary of an extremal set and hence its number of holes, a hole far from an end is excluded by comparing the Perron vector with a resolvent of the block, and a hole near an end by comparing with harmonic measures on a window next to it. The inequalities that enter are certified by computer.

The problem comes from the design of marine reserves along a coastline. If $C_{ij}$ is the fraction of larvae released at habitat patch $i$ that settle at patch $j$, and larvae that settle outside the reserves are lost, then the persistence of a reserve network $S$ is governed by the spectral radius of the principal submatrix $C[S]$~\cite{HB2006}, and the usual perturbative choice of $S$ carries no optimality guarantee~\cite{JJ2011}. For a symmetric dispersal kernel of finite range on a uniform coastline the problem is exactly $\rho_m(\omega)$, which is why we allow weights, and Theorem~\ref{thm:all}(1) says that once the number of reserves is large, a contiguous network is the unique optimal choice for every such kernel whose support has greatest common divisor $1$. We return to this in Section~\ref{sec:remarks}.

The paper is organized as follows. Section~\ref{sec:def} collects the random walk tools and the harmonic functions of admissible sets. Section~\ref{sec:main} contains the proofs: Section~\ref{sec:half} proves Theorem~\ref{thm:main}, Section~\ref{sec:ends} builds approximate eigenvectors of sets whose holes lie near the ends and proves Proposition~\ref{prop:ends}, Section~\ref{sec:holes} excludes holes away from the ends and completes the proof of Theorem~\ref{thm:all}(1), and Section~\ref{sec:exact} proves Theorem~\ref{thm:all}(2). Section~\ref{sec:remarks} contains concluding remarks.

\section{Definitions and Tools}\label{sec:def}

Throughout Section~\ref{sec:def} and Sections~\ref{sec:half}--\ref{sec:holes}, $D$ is a finite set of positive integers with $|D|\ge2$ and $\gcd(D)=1$, $b=\max D$, and $\omega\colon D\to(0,\infty)$ is a weight, extended to $\Z$ as in the introduction; we put $W=2\sum_{k\in D}\omega_k$, $\omega_{\min}=\min_D\omega$ and $\omega_{\max}=\max_D\omega$. Let $(X_t)_{t\ge0}$ be the random walk on $\Z$ with transition probabilities $p(x,y)=\omega_{x-y}/W$, let $\sigma^2=\sum_{k\in D}\omega_kk^2\big/\sum_{k\in D}\omega_k$ be the variance of a step, and put $\kappa=2/\sigma^2$. The walk is symmetric, has mean zero and bounded steps, and is therefore recurrent; since $\gcd(D)=1$ it is irreducible on $\Z$. For the standard facts about random walks used below we refer to~\cite{LL} and~\cite{Spitzer}.

For a finitely supported $x\colon\Z\to\mathbb{R}$,
\[
W\|x\|^2-\sum_{i,j}A_{ij}x_ix_j=\sum_{k\in D}\omega_k\sum_{i\in\Z}(x_i-x_{i+k})^2,
\]
where $A$ is the adjacency matrix of $G(\Z,\omega)$, because $\sum_i(x_i-x_{i+k})^2=2\|x\|^2-2\sum_ix_ix_{i+k}$. Hence
\[
W-\rho(S)=\min\Big\{\sum_{k\in D}\omega_k\sum_{i\in\Z}(x_i-x_{i+k})^2:\ \operatorname{supp}x\subseteq S,\ \|x\|=1\Big\}.
\]
In Fourier terms the right side is $\frac{1}{2\pi}\int_{-\pi}^{\pi}|\hat x(\theta)|^2\bigl(4\sum_{k\in D}\omega_k\sin^2\frac{k\theta}{2}\bigr)\,d\theta$, and the weight vanishes on $[-\pi,\pi]$ only at $\theta=0$ because $\gcd(D)=1$. When $D=\{d\}$, it vanishes at every multiple of $2\pi/d$, which is why arithmetic progressions are extremal in that case.

\begin{defn}
A set $S\subseteq\Z$ is \emph{admissible} if it is bounded below and contains $[N,\infty)$ for some $N$. We write $T=\Z\setminus S$ and $\tau_S=\min\{t\ge0: X_t\in T\}$ for the exit time of $S$. An admissible set is a \emph{half-line} if $S=[s,\infty)$ for some $s$.
\end{defn}

By recurrence $\tau_S<\infty$ almost surely. Since every step has length at most $b$ and $T\supseteq(-\infty,\min S-1]$, we have $\min S-b\le X_{\tau_S}\le\max T$. Define
\[
h_S(i)=i-\E_i[X_{\tau_S}]\qquad(i\in\Z).
\]

\begin{lem}\label{lem:h}
Let $S$ be admissible.
\begin{enumerate}
\item $h_S=0$ on $T$, and $h_S(i)=\sum_jp(i,j)h_S(j)$ for every $i\in S$.
\item $h_S(i)-i$ is bounded, and $h_S\ge0$.
\item If $i\in S$ is joined to $[N,\infty)$ by a path in the subgraph of $G(\Z,D)$ induced by $S$, then $h_S(i)>0$.
\item $h_S$ is the only function with properties (1) and the boundedness of $h_S(i)-i$.
\item The limit $\beta(S)=\lim_{i\to\infty}(h_S(i)-i)$ exists.
\item For every integer $t$, $h_{S-t}(i)=h_S(i+t)$; in particular $\beta(S-t)=\beta(S)+t$.
\end{enumerate}
\end{lem}

\begin{proof}
(1) If $i\in T$ then $\tau_S=0$ and $h_S(i)=0$. If $i\in S$ then $\tau_S\ge1$, so by the Markov property $\E_i[X_{\tau_S}]=\sum_jp(i,j)\E_j[X_{\tau_S}]$, and $i=\sum_jp(i,j)j$ because the walk has mean zero.

(2) Boundedness follows from the bounds on $X_{\tau_S}$. Since $h_S(i)\to\infty$ and $S$ is bounded below, if $h_S$ took a negative value on $S$ it would attain a negative minimum at some $i_0\in S$. By (1), every $j$ with $|j-i_0|\in D$ would satisfy $h_S(j)=h_S(i_0)<0$, so $j\in S$. Repeating, every vertex reachable from $i_0$ by steps in $\pm D$ would lie in $S$; by irreducibility this includes points of $T$, a contradiction.

(3) If $h_S(i)=0$ for some $i\in S$, then by (1) and (2) every neighbor of $i$ in $S$ also has value $0$. Following the path to $[N,\infty)$ and then the steps $+\min D$, we would obtain $h_S=0$ at arbitrarily large integers, contradicting $h_S(i)\to\infty$.

(4) The difference $\chi$ of two such functions is bounded, harmonic on $S$ and zero on $T$. Then $\chi(X_{t\wedge\tau_S})$ is a bounded martingale, and optional stopping gives $\chi(i)=\E_i[\chi(X_{\tau_S})]=0$.

(5) This is Lemma~\ref{lem:renewal} applied to $i\mapsto \E_i[X_{\tau_S}]$.

(6) The walk started at $i$ and stopped on leaving $S-t$ is the walk started at $i+t$ and stopped on leaving $S$, shifted by $-t$.
\end{proof}

The next lemma is a quantitative form of the renewal theorem for the ladder heights of the walk~\cite[Chapter~XIII]{Feller}; we include the short proof because the geometric rate is used later.

\begin{lem}\label{lem:renewal}
Let $R$ be admissible and let $\chi$ be bounded and harmonic on $R$. Then $L=\lim_{i\to\infty}\chi(i)$ exists, and there are $C>0$ and $\varrho\in(0,1)$, where $\varrho$ depends only on $\omega$, such that $|\chi(i)-L|\le C\varrho^{\,i}$ for all $i\ge0$.
\end{lem}

\begin{proof}
Let $K=\max(\Z\setminus R)$ and $I=(K,K+b]$, so $I\subseteq R$. For $i>K+b$, let $\zeta=\min\{t\ge0: X_t\le K+b\}$. Since the walk is recurrent and irreducible, $\liminf_tX_t=-\infty$ almost surely, so $\zeta<\infty$ almost surely. Because $X_{\zeta-1}\ge K+b+1$ and every step is at least $-b$, we have $X_\zeta\in I$. For $t<\zeta$ the walk stays in $(K+b,\infty)\subseteq R$, where $\chi$ is harmonic, so $\chi(X_{t\wedge\zeta})$ is a bounded martingale and optional stopping gives
\begin{equation}\label{eq:window}
\chi(i)=\sum_{r=0}^{b-1}\Prob_i(X_\zeta=K+b-r)\,\chi(K+b-r).
\end{equation}
By translation invariance, $\Prob_i(X_\zeta=K+b-r)=\Prob(U_n=r)$ with $n=i-K-b$, where $U_n$ is defined as follows. Start the walk at $0$, and let $0=\Lambda_0<\Lambda_1<\Lambda_2<\cdots$ be the successive strict descending ladder heights, that is, $-\Lambda_j$ is the $j$-th strict new minimum of the walk. For the walk started at $n\ge1$, the first time it reaches $(-\infty,0]$ is the first ladder epoch with $\Lambda_j\ge n$, and its position there is $n-\Lambda_j$. Thus $U_n=\Lambda_J-n$ with $J=\min\{j:\Lambda_j\ge n\}$, and it suffices to show that the law of $U_n$ converges geometrically fast.

By the strong Markov property at the ladder epochs, which are almost surely finite, the increments $H_j=\Lambda_j-\Lambda_{j-1}$ are independent with a common law $f_k=\Prob(H_1=k)$. Every step has length at most $b$, so $f_k=0$ unless $1\le k\le b$. A first step of $-k$ with $k\in D$ gives $H_1=k$, so $f_k\ge\omega_k/W>0$ for every $k\in D$, and hence $\gcd\{k:f_k>0\}=1$. Let $v_j=\Prob(\Lambda_i=j\text{ for some }i)$, so $v_0=1$, and let $\nu=\sum_kkf_k$. If $U_n=r$, then the last ladder height below $n$ is some $j\le n-1$ and the next increment is $n+r-j\le b$. Hence, for $n\ge b$,
\begin{equation}\label{eq:overshoot}
\Prob(U_n=r)=\sum_{s=r+1}^{b}f_s\,v_{n+r-s}\qquad(0\le r\le b-1).
\end{equation}

The renewal sequence has generating function $V(z)=\sum_{j\ge0}v_jz^j=1/(1-F(z))$ for $|z|<1$, where $F(z)=\sum_{k=1}^bf_kz^k$ is a polynomial with $F(1)=1$. We locate the zeros of $1-F$ in $|z|\le1$. If $|z|\le1$ and $F(z)=1$, then $1=|F(z)|\le\sum_kf_k|z|^k\le1$, so $|z|=1$ and every $z^k$ with $f_k>0$ is a nonnegative real number of modulus one, that is, $z^k=1$. Since these $k$ have greatest common divisor $1$, we get $z=1$. Moreover $F'(1)=\nu>0$, so $z=1$ is a simple zero. Let $r_0>1$ be the smallest modulus of the other zeros of the polynomial $1-F$, and put $r_0=\infty$ if there are none. Near $z=1$ we have $1-F(z)=\nu(1-z)+O((1-z)^2)$, so
\[
Q(z)=\frac{1}{1-F(z)}-\frac{1}{\nu(1-z)}
\]
has a removable singularity at $z=1$ and is analytic in $|z|<r_0$. Fix $\varrho\in(1/r_0,1)$; it depends only on $\omega$. By the Cauchy estimates on $|z|=1/\varrho$, the coefficients $q_j$ of $Q$ satisfy $|q_j|\le C_0\varrho^j$. Comparing coefficients, $v_j=\frac{1}{\nu}+q_j$, so $|v_j-\frac{1}{\nu}|\le C_0\varrho^j$.

Put $\pi(r)=\frac{1}{\nu}\sum_{s=r+1}^bf_s$. By \eqref{eq:overshoot}, for $n\ge b$,
\[
\bigl|\Prob(U_n=r)-\pi(r)\bigr|\le\sum_{s=r+1}^bf_s\,C_0\varrho^{\,n+r-s}\le C_0\varrho^{\,n-b}.
\]
Let $L=\sum_{r=0}^{b-1}\pi(r)\,\chi(K+b-r)$ and $\|\chi\|=\sup|\chi|$. By \eqref{eq:window}, for $i\ge K+2b$,
\[
|\chi(i)-L|\le\sum_{r=0}^{b-1}\bigl|\Prob(U_{i-K-b}=r)-\pi(r)\bigr|\,\|\chi\|\le b\,C_0\,\|\chi\|\,\varrho^{\,i-K-2b}.
\]
Enlarging the constant to cover the finitely many $0\le i<K+2b$, where $|\chi(i)-L|\le2\|\chi\|$, gives the claim.
\end{proof}

For admissible $X$ and $x,y\in\Z$, let $G_X(x,y)$ be the expected number of visits to $y$ at times $0\le t<\tau_X$ for the walk started at $x$; this is the \emph{Green function} of the walk killed on leaving $X$. For $x\in X$, $G_X(x,x)=1/\Prob_x(\text{the walk leaves }X\text{ before returning to }x)$, which is finite. The function $G_X$ is symmetric, invariant under translating $X$, $x$ and $y$ together, and monotone: $X_1\subseteq X_2$ implies $G_{X_1}\le G_{X_2}$.

\begin{lem}\label{lem:monotone}
If $S_1\subseteq S_2$ are admissible, then $h_{S_1}\le h_{S_2}$ pointwise.
\end{lem}

\begin{proof}
The function $\chi=h_{S_2}-h_{S_1}$ is bounded by Lemma~\ref{lem:h}(2), harmonic on $S_1$, and equal to $h_{S_2}\ge0$ on $\Z\setminus S_1$. Optional stopping at $\tau_{S_1}$ gives $\chi(i)=\E_i[h_{S_2}(X_{\tau_{S_1}})]\ge0$.
\end{proof}

For admissible $R$ and $y\in\Z$ put
\[
\varphi_R(y)=\sum_{x\in R}p(y,x)\,h_R(x)=\frac{1}{W}\sum_{x\in R}\omega_{x-y}\,h_R(x),
\]
the \emph{harmonic measure from infinity} of $R$. Only finitely many $y\notin R$ have $\varphi_R(y)\ne0$.

\begin{lem}\label{lem:hm}
Let $R$ be admissible and let $\chi$ be bounded and harmonic on $R$. Then
\[
\lim_{i\to\infty}\chi(i)=\kappa\sum_{y\notin R}\varphi_R(y)\,\chi(y).
\]
\end{lem}

\begin{proof}
Write $\Delta f(x)=\sum_yp(x,y)(f(y)-f(x))$ and fix a large integer $M$. Put $B=R\cap(-\infty,M]$. Since $\Delta\chi=\Delta h_R=0$ on $B$,
\[
0=\sum_{x\in B}\bigl(\chi(x)\Delta h_R(x)-h_R(x)\Delta\chi(x)\bigr)=\sum_{x\in B}\sum_{y\notin B}p(x,y)\bigl(\chi(x)h_R(y)-h_R(x)\chi(y)\bigr),
\]
because the terms with $x,y\in B$ cancel in pairs by the symmetry of $p$, and the diagonal terms cancel. The pairs with $y\notin R$ contribute $-\sum_{y\notin R}\varphi_R(y)\chi(y)$, since $h_R(y)=0$ there. The remaining pairs have $x\le M<y$ and $y-x\in D$. Write $\chi=L+\xi$ and $h_R(z)=z+\beta(R)+\xi'(z)$, where $\xi$ and $\xi'$ decay geometrically by Lemma~\ref{lem:renewal}. Then
\[
\chi(x)h_R(y)-h_R(x)\chi(y)=L\,(y-x)+O\bigl(M(|\xi(x)|+|\xi(y)|)+|\xi'(x)|+|\xi'(y)|\bigr),
\]
and the error terms vanish as $M\to\infty$. For each $k\in D$ there are exactly $k$ pairs with $x\le M<y$ and $y-x=k$, each with $p(x,y)=\omega_k/W$, so these pairs contribute $L\sum_{k\in D}\frac{\omega_kk^2}{W}+o(1)=\frac{\sigma^2}{2}L+o(1)$. Letting $M\to\infty$ gives the claim.
\end{proof}

Finally, for finite $S\subseteq\Z$ we write $A_S$ for the adjacency matrix of the subgraph induced by $S$, and $A_{S,T}$ for the submatrix of $A$ with rows in $S$ and columns in $T$. A \emph{hole} of $S$ is a point $w\notin S$ with $\min S<w<\max S$, and its \emph{depth} is $\min(w-\min S,\max S-w)$.

\section{Main Results}\label{sec:main}

In this section, we prove Theorems~\ref{thm:main} and~\ref{thm:all}. The four subsections follow the outline given in the introduction.

\subsection{The half-line comparison}\label{sec:half}

By Lemma~\ref{lem:h}(6), $\beta([N-c,\infty))=\beta([0,\infty))+c-N$, so Theorem~\ref{thm:main} says that among all ends with the same number of vertices, the contiguous end pushes the boundary farthest. The proof moves one vertex at a time.

\begin{proof}[Proof of Theorem~\ref{thm:main}]
We first prove the inequality and its equality case. Suppose $S$ is not a half-line. Let $v=\min S$ and let $w=\max(\Z\setminus S)$; then $v<w$, every vertex above $w$ lies in $S$, and every vertex below $v$ lies outside $S$. Put
\[
R=S\setminus\{v\},\qquad S^*=R\cup\{w\}.
\]
The set $S^*$ is admissible and has the same $c$ vertices below $N$ as $S$ except that $v$ is replaced by $w$; since $v<w<N$, the sum of the elements of $S\cap(-\infty,N)$ strictly increases. This sum is at most the sum of $[N-c,N)$, so finitely many moves turn $S$ into $[N-c,\infty)$. It therefore suffices to show $\beta(S^*)>\beta(S)$.

\emph{Claim 1.} $\beta(S)-\beta(R)=\kappa\,\varphi_R(v)^2\,G_S(v,v)$ and $\beta(S^*)-\beta(R)=\kappa\,\varphi_R(w)^2\,G_{S^*}(w,w)$.

We prove the first identity; the second is identical with $(S,v)$ replaced by $(S^*,w)$. Let $\chi=h_S-h_R$. It is bounded, harmonic on $R$, equal to $h_S(v)$ at $v$, and zero on $\Z\setminus S$. Lemma~\ref{lem:hm} gives $\beta(S)-\beta(R)=\lim\chi=\kappa\,\varphi_R(v)h_S(v)$. To compute $h_S(v)$, optional stopping at $\tau_R$ gives $\chi(y)=h_S(v)\Prob_y(X_{\tau_R}=v)$ for $y\in R$. Since $h_S$ is harmonic at $v$,
\[
h_S(v)=\sum_{y\in R}p(v,y)\bigl(h_R(y)+\chi(y)\bigr)=\varphi_R(v)+h_S(v)\,\Prob_v(\text{return to }v\text{ before leaving }S),
\]
and therefore $h_S(v)=\varphi_R(v)\,G_S(v,v)$.

\emph{Claim 2.} $G_{S^*}(w,w)>G_S(v,v)$.

Since $S\subseteq[v,\infty)$ and $S^*\supseteq[w,\infty)$, monotonicity and translation invariance give
\[
G_S(v,v)\le G_{[v,\infty)}(v,v)=G_{[w,\infty)}(w,w)\le G_{S^*}(w,w).
\]
The first inequality is strict. Indeed, pick $z\in[v,\infty)\setminus S$. Because $\gcd(D)=1$ and the walk can move upward freely, there is a path of steps in $\pm D$ from $v$ to $z$ that stays in $[v+1,\infty)$ after its first step. Since the subgraph of $G(\Z,D)$ induced by $[v,\infty)$ is connected, there is also a path from $z$ back to $v$ inside $[v,\infty)$. With positive probability the walk follows the two paths in turn: it leaves $S$ before returning to $v$, but does not leave $[v,\infty)$ before returning to $v$. Since $S\subseteq[v,\infty)$, the escape probability from $v$ is therefore strictly larger for $S$ than for $[v,\infty)$.

\emph{Claim 3.} $\varphi_R(w)\ge\varphi_R(v)$, and $\varphi_R(w)>0$.

Put $t=w-v\ge1$. Since $R\subseteq[v+1,\infty)$ and $R-t\supseteq[w+1-t,\infty)=[v+1,\infty)$, we have $R\subseteq R-t$. Lemma~\ref{lem:monotone} and Lemma~\ref{lem:h}(6) give $h_R(i+t)=h_{R-t}(i)\ge h_R(i)$ for every $i$. All neighbors of $v$ in $R$ lie in $v+D$, while $w+D\subseteq[w+1,\infty)\subseteq R$. Hence
\[
\varphi_R(w)\ge\frac{1}{W}\sum_{k\in D}\omega_kh_R(w+k)\ge\frac{1}{W}\sum_{k\in D}\omega_kh_R(v+k)=\varphi_R(v).
\]
Each $w+k$ lies in $[w+1,\infty)$, which is joined to $[N,\infty)$ inside $R$, so $h_R(w+k)>0$ by Lemma~\ref{lem:h}(3), and $\varphi_R(w)>0$.

By Claims 1--3, if $\varphi_R(v)>0$ then $\beta(S^*)-\beta(R)>\beta(S)-\beta(R)$, and if $\varphi_R(v)=0$ then $\beta(S)=\beta(R)<\beta(S^*)$. In both cases $\beta(S^*)>\beta(S)$.

We now prove the quantitative bound. Write $\Delta(S)=\beta([N-c,\infty))-\beta(S)$, put $\beta_0=\beta([0,\infty))$, $R_j=[0,\infty)\setminus\{j\}$, and
\[
c(j)=\beta_0-\beta(R_j)\qquad(j\ge0).
\]
Applying Lemma~\ref{lem:hm} to $\chi=h_0-h_{R_j}$ and the argument of Claim~1 to the pair $(R_j\cup\{j\},j)$ gives $c(j)=\kappa\,h_0(j)^2/G_0(j,j)$; since $R_0=[1,\infty)$, Lemma~\ref{lem:h}(6) gives $c(0)=1$. Thus $\gamma(\omega)=c(1)-1$.

\emph{Claim 4.} $c(j)$ is strictly increasing in $j$.

Since $R_{j+1}-1=R_j\cup\{-1\}$, Lemma~\ref{lem:h}(6) gives $\beta(R_{j+1})=\beta(R_j\cup\{-1\})-1$, so it suffices to show $\beta(R_j\cup\{-1\})-\beta(R_j)<1$. By Claim~1,
\begin{align*}
\beta(R_j\cup\{-1\})-\beta(R_j)&=\kappa\,\varphi_{R_j}(-1)^2\,G_{R_j\cup\{-1\}}(-1,-1),\\
1&=\kappa\,\varphi_{[0,\infty)}(-1)^2\,G_{[-1,\infty)}(-1,-1),
\end{align*}
where the second identity is $\beta([-1,\infty))-\beta([0,\infty))=1$. Since $R_j\subsetneq[0,\infty)$, Lemma~\ref{lem:monotone} gives $\varphi_{R_j}(-1)\le\varphi_{[0,\infty)}(-1)$, and the argument of Claim~2 gives
\[
G_{R_j\cup\{-1\}}(-1,-1)<G_{[-1,\infty)}(-1,-1).
\]
Finally $\varphi_{R_j}(-1)>0$, because $|D|\ge2$ gives $-1$ a neighbor in $R_j\cap[0,\infty)$ other than $j$, and $h_{R_j}$ is positive there by Lemma~\ref{lem:h}(3).

\emph{Claim 5.} For $d\ge2$ put $U_d=\{-d\}\cup[0,\infty)$. Then $\Delta(U_d)$ is nondecreasing in $d$, and $\Delta(U_3)>\Delta(U_2)$.

By Claim~1, $\Delta(U_d)=1-\kappa\,\varphi_0(-d)^2\,G_{U_d}(-d,-d)$ with $\varphi_0=\varphi_{[0,\infty)}$. The function $h_0$ is strictly increasing: $h_0(i+1)-h_0(i)=h_{[-1,\infty)}(i)-h_0(i)>0$ by Lemma~\ref{lem:h}(6) and Lemma~\ref{lem:monotone}, with strictness from optional stopping at the exit from $[0,\infty)$. Hence $\varphi_0(-d)=\frac{1}{W}\sum_{k\in D,\,k\ge d}\omega_kh_0(k-d)$ is nonincreasing in $d$, and $\varphi_0(-3)<\varphi_0(-2)$: if $D$ has an element $k\ge3$, the term $h_0(k-3)<h_0(k-2)$ decreases strictly, and otherwise $\varphi_0(-3)=0<\varphi_0(-2)$, since $D$ has an element $k\ge2$ and $h_0(k-2)>0$. Since $U_{d+1}+1=\{-d\}\cup[1,\infty)\subseteq U_d$, translation invariance and monotonicity give $G_{U_{d+1}}(-d-1,-d-1)\le G_{U_d}(-d,-d)$.

By the first part of the proof, $\Delta(S)$ is the sum of the increments $\beta(S^*)-\beta(S)$ over the moves that turn $S$ into $[N-c,\infty)$, and each increment is positive. The set before the last move has the form, after translation,
\[
S_{v,w}=\{v\}\cup\bigl([0,\infty)\setminus\{w\}\bigr),\qquad v\le-1,\ w\ge0,
\]
and the last increment equals $\Delta(S_{v,w})=\beta_0-\beta(S_{v,w})$. Hence $\Delta(S)\ge\Delta(S_{v,w})$, with equality only if a single move is made, that is, only if $S$ itself is a translate of some $S_{v,w}$. We show $\Delta(S_{v,w})\ge c(1)-1$, with equality only for $(v,w)=(-1,0)$.

If $v=-1$, then $S_{-1,w}=R_{w+1}-1$, so $\beta(S_{-1,w})=\beta_0-c(w+1)+1$ and $\Delta(S_{-1,w})=c(w+1)-1\ge c(1)-1$ by Claim~4, with equality only for $w=0$.

If $v\le-2$ and $w\ge1$, put $S^\sharp=S_{v,w}\cup\{w\}=\{v\}\cup[0,\infty)$. By Claim~1 applied to $(S^\sharp,w)$,
\[
\beta(S_{v,w})=\beta(S^\sharp)-\kappa\,\varphi_{S_{v,w}}(w)^2\,G_{S^\sharp}(w,w).
\]
Since $S_{v,w}\supseteq R_w$ and $S^\sharp\supseteq[0,\infty)$, monotonicity gives
\[
\varphi_{S_{v,w}}(w)^2\,G_{S^\sharp}(w,w)\ge\varphi_{R_w}(w)^2\,G_0(w,w)=h_0(w)^2/G_0(w,w),
\]
so $\beta(S_{v,w})\le\beta(S^\sharp)-c(w)$. The first part of the theorem gives $\beta(S^\sharp)<\beta([-1,\infty))=\beta_0+1$, hence $\Delta(S_{v,w})>c(w)-1\ge c(1)-1$.

If $v\le-2$ and $w=0$, then $S_{v,0}=U_{1-v}+1$ with $1-v\ge3$, so $\Delta(S_{v,0})=\Delta(U_{1-v})>\Delta(U_2)$ by Claim~5. Since $U_2+1=S_{-1,0}$, we have $\Delta(U_2)=\Delta(S_{-1,0})=c(1)-1$.

Positivity of $\gamma(\omega)=c(1)-c(0)$ is Claim~4.
\end{proof}

\begin{rem}\label{rem:gamma}
The proof uses only $|D|\ge2$ and $\gcd(D)=1$, through irreducibility in Claim~2 and Lemma~\ref{lem:h}, and neither the structure of $D$, nor the values of the weights, nor a closed form of $\beta([0,\infty))$. The extremal end $\{-1\}\cup[1,\infty)$ has a single vertex separated from the half-line by one hole: removing the second vertex of a contiguous end costs $c(1)$, whereas shortening the end by one vertex costs only $c(0)=1$. In the unweighted case, numerically, $\gamma\approx0.4597$, $0.3938$ and $0.3238$ for $D=\{2,3\}$, $\{2,3,4\}$ and $\{3,4,5\}$. The number $\gamma$ is the limit as $\vartheta\to0$ of $H(\vartheta)/\vartheta$, where $H$ is the phase gap in Proposition~\ref{prop:satellite} below, which compares the same two ends at positive frequency; the constants $H_{\mathrm{lo}}$ in Table~\ref{tab:sat} are certified lower bounds for $H(\vartheta)/\vartheta$ on $(0,\theta_{\max}]$.
\end{rem}

\subsection{End patterns}\label{sec:ends}

In this subsection we show that among the sets whose holes lie within a bounded distance of the ends, the block has the largest spectral radius once $m$ is large. The tool is an approximate eigenvector glued from generalized eigenfunctions of the two ends, whose phases are governed by the offsets of Theorem~\ref{thm:main}.

Let $f(\theta)=2\sum_{k\in D}\omega_k\cos k\theta$, so that $f(\theta)$ is the eigenvalue of $A$ on the plane wave $e^{i\theta j}$. Put $F(z)=\sum_{k\in D}\omega_k(z^k+z^{-k})$ and $Q_\vartheta(z)=z^b\bigl(F(z)-f(\vartheta)\bigr)$, a real polynomial of degree $2b$ with $Q_\vartheta(0)=\omega_b>0$ and $z^{2b}Q_\vartheta(1/z)=Q_\vartheta(z)$.

\begin{lem}\label{lem:roots}
There are $\vartheta_1\in(0,\pi/(2b)]$ and $r_1\in(0,1)$ such that for every $\vartheta\in[0,\vartheta_1]$:
\begin{enumerate}
\item the roots of $Q_\vartheta$ on the unit circle are $e^{\pm i\vartheta}$, simple if $\vartheta>0$, and a double root $1$ if $\vartheta=0$;
\item $Q_\vartheta$ has exactly $b-1$ roots, counted with multiplicity, in $|z|<1$, all of them in $|z|<r_1$.
\end{enumerate}
\end{lem}

\begin{proof}
For $z=e^{i\theta}$ we have $Q_\vartheta(z)=z^b(f(\theta)-f(\vartheta))$. If $f(\theta)=W$, then $\cos k\theta=1$ for every $k\in D$, since every $\omega_k>0$, so $k\theta\in2\pi\Z$ for every $k\in D$, and $\theta\in2\pi\Z$ because $\gcd(D)=1$. Hence $\epsilon_0=W-\max\{f(\theta):\pi/(2b)\le|\theta|\le\pi\}>0$. On $(0,\pi/(2b)]$ we have $f'(\theta)=-2\sum_k\omega_kk\sin k\theta<0$. Choose $\vartheta_1\le\pi/(2b)$ with $f(\vartheta_1)>W-\epsilon_0$. For $\vartheta\le\vartheta_1$ and $\theta\in(-\pi,\pi]$, $f(\theta)=f(\vartheta)$ forces $|\theta|<\pi/(2b)$ and then $\theta=\pm\vartheta$. If $\vartheta>0$, then $\frac{d}{d\theta}F(e^{i\theta})=iz F'(z)=f'(\theta)\ne0$ at $\theta=\vartheta$, so $e^{\pm i\vartheta}$ are simple roots. If $\vartheta=0$, then $F(z)-W=\sum_k\omega_kz^{-k}(z^k-1)^2$, which vanishes to order exactly two at $z=1$, since $(z-1)^{-2}(F(z)-W)\to\sum_k\omega_kk^2$.

The other $2b-2$ roots lie off the unit circle and are closed under $z\mapsto1/z$, so $b-1$ of them lie inside. At $\vartheta=0$ these have modulus at most some $\bar r<1$. The roots of $Q_\vartheta$ depend continuously on $\vartheta$ as a multiset, and the two roots on the circle stay there, so after shrinking $\vartheta_1$ the inside roots stay in $|z|<r_1$ for any fixed $r_1\in(\bar r,1)$.
\end{proof}

For $\vartheta\in[0,\vartheta_1]$ let $E_\vartheta$ be the space of real sequences $y\colon\Z\to\mathbb R$ with $(Ay)(i)=f(\vartheta)y(i)$ for every $i\in\Z$ and $y(j)\to0$ as $j\to\infty$. A solution of this recurrence of order $2b$ is determined by its values on any $2b$ consecutive integers, and it is a sum of terms $p(j)z^j$ over the roots $z$ of $Q_\vartheta$, with $\deg p$ less than the multiplicity. Hence $\dim E_\vartheta=b-1$, and $E_\vartheta$ consists of the solutions built from the roots in $|z|<r_1$.

Let $M_\vartheta$ be the companion matrix of the recurrence; it maps the window $(y(j),\dots,y(j+2b-1))$ of a solution to the next window. Fix $r_2$ with $r_1<r_2<1$. The Riesz projection
\[
\mathsf P_\vartheta=\frac1{2\pi i}\oint_{|z|=r_2}(zI-M_\vartheta)^{-1}dz
\]
onto the window space of $E_\vartheta$ is real and depends analytically on $\vartheta$ near $[0,\vartheta_1]$, because $M_\vartheta$ depends analytically on $f(\vartheta)$ and has no eigenvalue on $|z|=r_2$. Choose real vectors $v_1,\dots,v_{b-1}$ such that $\mathsf P_0v_l$ are linearly independent. After shrinking $\vartheta_1$, the sequences $e_l(\vartheta)\in E_\vartheta$ whose window $(y(-b),\dots,y(b-1))$ is $\mathsf P_\vartheta v_l$ form a basis of $E_\vartheta$ that is analytic in $\vartheta$, and there is $C_E$ with
\begin{equation}\label{eq:decay}
|e_l(\vartheta)(j)|\le C_E\,r_2^{\,j}\qquad(j\ge-b,\ 1\le l\le b-1,\ 0\le\vartheta\le\vartheta_1),
\end{equation}
since $M_\vartheta^n\mathsf P_\vartheta=\frac1{2\pi i}\oint_{|z|=r_2}z^n(zI-M_\vartheta)^{-1}dz$ and the resolvent is bounded on $|z|=r_2$ uniformly for $\vartheta$ in the compact interval, so $\|M_\vartheta^n\mathsf P_\vartheta\|\le Cr_2^{\,n}$.

For a finite set $P\subseteq(-\infty,-1]$ put $S_P=P\cup[0,\infty)$; we call $P$ an \emph{end pattern}. Write $s_\vartheta(j)=\sin(\vartheta j)/\vartheta$, read as $j$ when $\vartheta=0$, and $c_\vartheta(j)=\cos\vartheta j$; both are entire functions of $\vartheta$ and solve $(Ay)(i)=f(\vartheta)y(i)$ on $\Z$.

\begin{lem}\label{lem:psi}
For every end pattern $P$ there is $\vartheta_P\in(0,\vartheta_1]$ such that for every $\vartheta\in[0,\vartheta_P]$ there is a unique function $\psi=\psi_{P,\vartheta}\colon\Z\to\mathbb R$ with
\begin{enumerate}
\item $\psi=0$ off $S_P$;
\item $(A\psi)(i)=f(\vartheta)\psi(i)$ for every $i\in S_P$;
\item $\psi(j)=s_\vartheta(j)+g\,c_\vartheta(j)+\eta(j)$ for every $j\ge-b$, for some $g\in\mathbb R$ and $\eta\in E_\vartheta$.
\end{enumerate}
The number $g=g_P(\vartheta)$, the coefficients of $\eta=\eta_{P,\vartheta}$ in the basis $e_l(\vartheta)$, and the values $\psi(j)$ for each fixed $j$ are analytic in $\vartheta\in[0,\vartheta_P]$. At $\vartheta=0$ we have $\psi=h_{S_P}$ and $g_P(0)=\beta(S_P)$.
\end{lem}

\begin{proof}
Write $P_<=P\cap(-\infty,-b-1]$. The unknowns are $g$, the coefficients $a_1,\dots,a_{b-1}$ of $\eta=\sum_la_le_l(\vartheta)$, and the values $\psi(p)$ for $p\in P_<$; on $[-b,\infty)$ the function $\psi$ is given by~(3). The conditions are $\psi(j)=0$ for $j\in[-b,-1]\setminus P$, and the eigenvalue equation~(2) at every $i\in P$. The eigenvalue equation at $i\ge0$ holds automatically: all neighbors of $i$ lie in $[-b,\infty)$, where $\psi$ agrees with a solution of the recurrence on $\Z$. There are $1+(b-1)+|P_<|$ unknowns and $(b-|P\cap[-b,-1]|)+|P|=b+|P_<|$ equations, so we have a square linear system whose coefficients are analytic in $\vartheta$; the term $s_\vartheta$ contributes the right-hand side.

At $\vartheta=0$ the system is nonsingular. Indeed, a solution of the homogeneous system gives a function $\psi_0$ that vanishes off $S_P$, is harmonic on $S_P$, and equals $g+\eta(j)$ on $[-b,\infty)$, so it is bounded. Optional stopping at the exit time from $S_P$, which is finite almost surely by recurrence, gives $\psi_0=0$. Then $\eta=-g$ on $[-b,\infty)$; since $\eta(j)\to0$, $g=0$ and $\eta$ vanishes on $2b$ consecutive integers, so $\eta=0$ and all $a_l=0$. By continuity the system is nonsingular on some $[0,\vartheta_P]$, and its solution is analytic there. The representation in~(3) is unique, because $c_\vartheta$ does not tend to $0$ and so is not in $E_\vartheta$; hence a function $\psi$ with (1)--(3) determines the unknowns, and uniqueness of $\psi$ follows from nonsingularity.

At $\vartheta=0$ the function $\psi$ vanishes off $S_P$, is harmonic on $S_P$, and $\psi(j)-j=g+\eta(j)$ is bounded on $[-b,\infty)$, hence on $S_P$, which differs from $[-b,\infty)$ by finitely many points. Since $S_P$ is admissible, Lemma~\ref{lem:h}(4) gives $\psi=h_{S_P}$, and $g=\lim_j(\psi(j)-j)=\beta(S_P)$.
\end{proof}

For $0<\vartheta\le\vartheta_P$ write
\[
s_\vartheta(j)+g_P(\vartheta)c_\vartheta(j)=A_P\sin(\vartheta j+\Phi_P),\qquad A_P=\bigl(\vartheta^{-2}+g_P^2\bigr)^{1/2},\quad\Phi_P=\arctan\bigl(\vartheta g_P\bigr),
\]
with $\Phi_P\in(-\pi/2,\pi/2)$; this holds because $A_P\cos\Phi_P=1/\vartheta$ and $A_P\sin\Phi_P=g_P$. Define the \emph{phase deficit}
\[
H_P(\vartheta)=\vartheta|P|+\Phi_\emptyset(\vartheta)-\Phi_P(\vartheta)\qquad(0<\vartheta\le\min(\vartheta_P,\vartheta_\emptyset)),
\]
and call $P$ \emph{contiguous} if $P=[-c,-1]$ for some $c\ge0$.

\begin{lem}\label{lem:deficit}
Let $P$ be an end pattern.
\begin{enumerate}
\item $H_P(\vartheta)/\vartheta\to\beta_0+|P|-\beta(S_P)$ as $\vartheta\to0$, where $\beta_0=\beta([0,\infty))$. If $P$ is not contiguous, this limit is at least $\gamma(\omega)$.
\item If $P$ is contiguous, then $H_P(\vartheta)=0$ for all sufficiently small $\vartheta>0$.
\end{enumerate}
\end{lem}

\begin{proof}
(1) $\Phi_P(\vartheta)/\vartheta=\arctan(\vartheta g_P(\vartheta))/\vartheta\to g_P(0)=\beta(S_P)$ by Lemma~\ref{lem:psi}, and $\beta(S_\emptyset)=\beta_0$. If $P$ is not contiguous, then $S_P$ is not a half-line, and Theorem~\ref{thm:main} with $N=0$ and $c=|P|$ gives $\beta(S_P)\le\beta([-|P|,\infty))-\gamma(\omega)=\beta_0+|P|-\gamma(\omega)$, using Lemma~\ref{lem:h}(6).

(2) Let $P=[-c,-1]$, so $S_P=[-c,\infty)$, and put $\phi(j)=\psi_{\emptyset,\vartheta}(j+c)$. It vanishes off $S_P$ and satisfies the eigenvalue equation on $S_P$. For $j\ge-b$ we have $j+c\ge-b$, so
\[
\phi(j)=s_\vartheta(j+c)+g_\emptyset c_\vartheta(j+c)+\eta_\emptyset(j+c)=\varsigma\,s_\vartheta(j)+\varsigma'c_\vartheta(j)+\eta_\emptyset(j+c),
\]
with $\varsigma=\cos\vartheta c-\vartheta g_\emptyset\sin\vartheta c$, and $\eta_\emptyset(\cdot+c)\in E_\vartheta$. For small $\vartheta$ we have $\varsigma>0$, and by uniqueness in Lemma~\ref{lem:psi}, $\psi_{P,\vartheta}=\phi/\varsigma$. The oscillating part of $\phi$ is $A_\emptyset\sin(\vartheta j+\vartheta c+\Phi_\emptyset)$, so the oscillating part of $\psi_{P,\vartheta}$ is a positive multiple of $\sin(\vartheta j+\vartheta c+\Phi_\emptyset)$. Since both $\Phi_P$ and $\vartheta c+\Phi_\emptyset$ lie in $(-\pi/2,\pi/2)$ for small $\vartheta$, $\Phi_P=\vartheta c+\Phi_\emptyset$, that is, $H_P=0$.
\end{proof}

An interval with two end patterns has an approximate eigenvector glued from the eigenfunctions of its ends.

\begin{lem}\label{lem:approx}
Let $P_L,P_R$ be end patterns, let $e\ge2b$, and let $S=P_L\cup[0,e]\cup(e-P_R)$. Let $0<\vartheta\le\min(\vartheta_{P_L},\vartheta_{P_R})$ satisfy
\[
\vartheta e+\Phi_{P_L}(\vartheta)+\Phi_{P_R}(\vartheta)=\pi.
\]
Write $\psi_L=\psi_{P_L,\vartheta}$, $\psi_R=\psi_{P_R,\vartheta}$, $\eta_L,\eta_R$ for their decaying parts, $\varpi=A_{P_L}/A_{P_R}$, and define $y\colon\Z\to\mathbb R$ by
\[
y(j)=\psi_L(j)+\varpi\,\psi_R(e-j)-A_{P_L}\sin(\vartheta j+\Phi_{P_L}).
\]
Then
\begin{equation}\label{eq:tworeps}
y(j)=\psi_L(j)+\varpi\,\eta_R(e-j)\quad(j\le e+b),\qquad y(j)=\varpi\,\psi_R(e-j)+\eta_L(j)\quad(j\ge-b).
\end{equation}
For $i\in S$ put $\delta_i=(A_Sy)(i)-f(\vartheta)y(i)$. Then $\delta_i=0$ unless $i$ has a neighbor outside $S$, and
\[
|\delta_i|\le W\max(1,\varpi)\,\sup_{n\ge e+1}\max\bigl(|\eta_L(n)|,|\eta_R(n)|\bigr).
\]
Consequently, if $S$ induces a connected subgraph and $y>0$ on $S$, then $|\rho(S)-f(\vartheta)|\le\max_{i\in S}|\delta_i|/y(i)$.
\end{lem}

\begin{proof}
For $j\le e+b$ we have $e-j\ge-b$, so Lemma~\ref{lem:psi}(3) gives $\psi_R(e-j)=A_{P_R}\sin(\vartheta(e-j)+\Phi_{P_R})+\eta_R(e-j)$. By the quantization condition, $\vartheta(e-j)+\Phi_{P_R}=\pi-(\vartheta j+\Phi_{P_L})$, so $\varpi A_{P_R}\sin(\vartheta(e-j)+\Phi_{P_R})=A_{P_L}\sin(\vartheta j+\Phi_{P_L})$, which gives the first identity in~\eqref{eq:tworeps}. The second is symmetric, using $\psi_L(j)=A_{P_L}\sin(\vartheta j+\Phi_{P_L})+\eta_L(j)$ for $j\ge-b$.

Let $i\in S$ with $i\le e-b$. All neighbors $t$ of $i$ satisfy $t\le e$, and $S\cap(-\infty,e]=S_{P_L}\cap(-\infty,e]$, so $i\in S_{P_L}$, and a neighbor $t\notin S$ satisfies $t\notin S_{P_L}$, hence $\psi_L(t)=0$, and $t<0$, hence $e-t\ge e+1$. Using the first identity, $(A\psi_L)(i)=f(\vartheta)\psi_L(i)$ and $(A\eta_R)(e-i)=f(\vartheta)\eta_R(e-i)$, we get
\[
(A_Sy)(i)-f(\vartheta)y(i)=-\sum_{t\notin S}\omega_{t-i}\,y(t)=-\varpi\sum_{t\notin S}\omega_{t-i}\,\eta_R(e-t).
\]
For $i\in S$ with $i\ge b$ the symmetric argument with the second identity gives
\[
\delta_i=-\sum_{t\notin S}\omega_{t-i}\,\eta_L(t),
\]
where every $t\notin S$ with $\omega_{t-i}>0$ satisfies $t\ge e+1$. Since $e\ge2b$, every $i\in S$ falls under one of the two cases, and in either case $\sum_t\omega_{t-i}=W$. The last statement is the Collatz--Wielandt inequality
\[
\min_{i\in S}\frac{(A_Sy)(i)}{y(i)}\le\rho(S)\le\max_{i\in S}\frac{(A_Sy)(i)}{y(i)}
\]
for the irreducible matrix $A_S$ and the positive vector $y|_S$.
\end{proof}

For $D=\{2,3\}$, $\{2,3,4\}$, $\{3,4,5\}$ the roots of $Q_\vartheta$ inside the unit disk are simple for $0\le\vartheta\le\theta_{\max}:=\pi/290$, so the decaying part is $\eta_{P,\vartheta}(j)=\sum_{l=1}^{b-1}c_lz_l^{\,j+b}$, where $z_l$ are these roots. In the certified computations of Section~\ref{sec:exact} the linear system of Lemma~\ref{lem:psi} was solved in ball arithmetic on a cover of $[0,\theta_{\max}]$, which shows $\vartheta_P\ge\theta_{\max}$ for the patterns used there; we write $g_{\mathrm{lo}},g_{\mathrm{hi}}$ for enclosures of $g_\emptyset$ on $[0,\theta_{\max}]$, $r$ for an upper bound of $\max_l|z_l|$ and $C$ for an upper bound of $\sum_l|c_l|$, so that $\sup_{n\ge e+1}|\eta(n)|\le Cr^{\,e+1+b}$ in Lemma~\ref{lem:approx}.

\begin{prop}\label{prop:ends}
For every integer $J\ge0$ there is $M=M(\omega,J)$ with the following property. Let $m\ge M$ and let $S\subseteq\Z$ be a set with $|S|=m$ that induces a connected subgraph of $G(\Z,D)$ and all of whose holes have depth at most $J$. Then $\rho(S)\le\rho(B_m)$, with equality only if $S$ is a translate of $B_m$.
\end{prop}

\begin{proof}
Fix $J\ge0$. Let $\mathcal P_J$ be the set of end patterns $P\subseteq[-J-1,-1]$ with $-J-1\in P$ such that every point of $P$ is joined to $[0,\infty)$ by a path in the subgraph induced by $S_P$. It is finite, contains the contiguous pattern $[-J-1,-1]$, and we add $P=\emptyset$ to it.

\emph{Uniform constants.} By Lemmas~\ref{lem:psi} and~\ref{lem:deficit} and finiteness of $\mathcal P_J$, there is $\vartheta_*\in(0,\vartheta_1]$ such that $\vartheta_*\le\vartheta_P$ for all $P\in\mathcal P_J$ and the following hold for all $P\in\mathcal P_J$ and $0<\vartheta\le\vartheta_*$.
\begin{enumerate}
\item[(a)] $H_P(\vartheta)\ge\gamma\vartheta/2$ if $P$ is not contiguous, and $H_P(\vartheta)=0$ if it is, where $\gamma=\gamma(\omega)$.
\item[(b)] $|g_P(\vartheta)|\le g_*$ for a constant $g_*$, and $\vartheta g_*\le\pi/4$; hence $1\le\vartheta A_P\le2$ and $|\Phi_P|\le\vartheta g_*$.
\item[(c)] $|\eta_{P,\vartheta}(n)|\le C\,r_2^{\,n}$ for $n\ge-b$, by~\eqref{eq:decay} and boundedness of the coefficients $a_l$.
\item[(d)] $\Phi_\emptyset$ is Lipschitz on $[0,\vartheta_*]$ with some constant $\ell$, being analytic.
\end{enumerate}
Fix $J_0\ge2g_*+2b$ with $Cr_2^{\,j}\le j/(8\pi)$ for $j\ge J_0$. Every $j\in S_P\cap(-\infty,J_0]$ is joined to $[0,\infty)$ inside $S_P$, so $h_{S_P}(j)>0$ by Lemma~\ref{lem:h}(3). Since $\psi_{P,\vartheta}(j)\to h_{S_P}(j)$ as $\vartheta\to0$, after shrinking $\vartheta_*$ we may also assume:
\begin{enumerate}
\item[(e)] $\psi_{P,\vartheta}(j)\ge c_0$ for all $j\in S_P\cap(-\infty,J_0]$, for a constant $c_0>0$.
\end{enumerate}

\emph{Positivity.} Let $P_L,P_R\in\mathcal P_J$, $e\ge 2J_0$, and let $\vartheta\le\vartheta_*$ satisfy the quantization condition of Lemma~\ref{lem:approx}. We show that, for large $e$, $y>0$ on $S$ and $y\ge c_0/4$ at the points of $S$ within $J_0$ of an end. By~(b), $1/2\le\varpi\le2$. Consider $j\le e/2$, where $y(j)=\psi_L(j)+\varpi\eta_R(e-j)$; the case $j\ge e/2$ is symmetric, with $y(j)=\varpi\psi_R(e-j)+\eta_L(j)$, where the main term is $\varpi A_{P_R}\sin(\cdot)=A_{P_L}\sin(\cdot)$ with $A_{P_L}\ge1/\vartheta$, the decaying term is $\eta_L$, and near the right end $\varpi\psi_R\ge c_0/2$ by~(e) and $\varpi\ge1/2$. If $j\le J_0$, then $y(j)\ge c_0-2Cr_2^{\,e-J_0}\ge c_0/4$ for large $e$. If $J_0<j\le e/2$, put $\xi=\vartheta j+\Phi_{P_L}$. By~(b), $\xi\ge\vartheta(j-g_*)\ge\vartheta j/2$, and $\xi\le\vartheta e/2+\Phi_{P_L}=\pi/2+(\Phi_{P_L}-\Phi_{P_R})/2\le\pi/2+\vartheta g_*\le3\pi/4$, so $\sin\xi\ge\frac2\pi\min(\vartheta j/2,\pi/4)$. Since $A_{P_L}\ge1/\vartheta$,
\[
y(j)\ge\frac1\pi\min\Bigl(j,\frac{\pi}{2\vartheta}\Bigr)-Cr_2^{\,j}-2Cr_2^{\,e/2}.
\]
The quantization condition and~(b) give $\vartheta\le(\pi+2\vartheta g_*)/e\le 2\pi/e$, so $\pi/(2\vartheta)\ge e/4\ge j/2$, and the right side is at least $j/(2\pi)-j/(8\pi)-2Cr_2^{\,e/2}>0$ for large $e$.

\emph{Spectral radius of an interval with two end patterns.} Let $S$ be as in the statement, $L=\min S$, $R=\max S$. Every point of $[L+J+1,R-J-1]$ is in $S$, since a hole there would have depth larger than $J$. Translate so that $L+J+1=0$ and put $e=R-J-1-(L+J+1)$, $P_L=S\cap[-J-1,-1]$ and $P_R=e-\bigl(S\cap[e+1,e+J+1]\bigr)$. Then $S=P_L\cup[0,e]\cup(e-P_R)$, and $-J-1\in P_L$, $-J-1\in P_R$. Since $S$ is connected, every point of $P_L$ is joined to $[0,e]$ by a path in $S$, and up to its first visit to $[0,\infty)$ such a path stays in $S\cap(-\infty,-1]=P_L$; hence $P_L\in\mathcal P_J$, and likewise $P_R\in\mathcal P_J$. Since the hull of $S$ has at least $m$ points, $e\ge m-2J-3$, and $|P_L|+|P_R|+e+1=m$.

Put
\[
\Gamma_S(\vartheta)=\vartheta e+\Phi_{P_L}(\vartheta)+\Phi_{P_R}(\vartheta)=\Gamma_B(\vartheta)-H_{P_L}(\vartheta)-H_{P_R}(\vartheta),\qquad \Gamma_B(\vartheta)=\vartheta(m-1)+2\Phi_\emptyset(\vartheta),
\]
where the identity uses $e=m-1-|P_L|-|P_R|$. The function $\Gamma_B$ is the quantization function of $B_m=[0,m-1]$ with empty patterns. Both are continuous on $[0,\vartheta_*]$, vanish at $0$, and satisfy $\Gamma_S(\vartheta)\ge\vartheta(e-2g_*)$ and $\Gamma_B(\vartheta)\ge\vartheta(m-1-2g_*)$ by~(b). For $e>2g_*+\pi/\vartheta_*$ let $\vartheta_S$ be the smallest root of $\Gamma_S=\pi$ in $(0,\vartheta_*]$; it exists by the intermediate value theorem. For $m-1>2\ell$ the function $\Gamma_B$ is strictly increasing by~(d), so $\Gamma_B=\pi$ has a unique root $\vartheta_B$, and $\pi/(m-1+2g_*)\le\vartheta_B\le\pi/(m-1-2g_*)$.

By the positivity argument and Lemma~\ref{lem:approx}, with~(b) and~(c),
\[
\rho(S)\le f(\vartheta_S)+\frac{8WC}{c_0}r_2^{\,e},\qquad\rho(B_m)\ge f(\vartheta_B)-\frac{8WC}{c_0}r_2^{\,m-1},
\]
because $\delta_i\ne0$ only for $i$ with a neighbor outside $S$, and such $i$ lie within $b\le J_0$ of an end, where $y\ge c_0/4$.

\emph{Comparison.} If $S$ is not a translate of $B_m$, then $P_L$ or $P_R$ is not contiguous, so by~(a), $H_{P_L}+H_{P_R}\ge\gamma\vartheta/2$ on $(0,\vartheta_*]$. Hence $\Gamma_B(\vartheta_S)=\pi+H_{P_L}(\vartheta_S)+H_{P_R}(\vartheta_S)\ge\pi+\gamma\vartheta_S/2>\pi$, so $\vartheta_S>\vartheta_B$, and by~(d)
\[
\frac{\gamma\vartheta_B}{2}\le\frac{\gamma\vartheta_S}2\le \Gamma_B(\vartheta_S)-\Gamma_B(\vartheta_B)\le(m-1+2\ell)(\vartheta_S-\vartheta_B).
\]
On $[0,\pi/(2b)]$ we have $-f'(\theta)=2\sum_k\omega_kk\sin k\theta\ge\frac4\pi\sum_k\omega_kk^2\,\theta$. Therefore
\[
f(\vartheta_B)-f(\vartheta_S)\ge\frac4\pi\sum_{k\in D}\omega_kk^2\,\vartheta_B(\vartheta_S-\vartheta_B)\ge\frac{2\gamma\sum_k\omega_kk^2}{\pi}\cdot\frac{\vartheta_B^2}{m-1+2\ell}\ge\frac{c_1}{m^3}
\]
for a constant $c_1>0$ depending only on $\omega$ and $J$. Combining,
\[
\rho(B_m)-\rho(S)\ge\frac{c_1}{m^3}-\frac{8WC}{c_0}\bigl(r_2^{\,m-2J-3}+r_2^{\,m-1}\bigr)>0
\]
for all $m$ beyond a bound $M(\omega,J)$, which also absorbs the finitely many conditions on $e\ge m-2J-3$ used above. If $S$ is a translate of $B_m$, then $\rho(S)=\rho(B_m)$.
\end{proof}

\begin{rem}
The proof suggests that the worst end pattern is the satellite $\{-2\}\cup[0,\infty)$, which attains $\gamma(\omega)$ in Theorem~\ref{thm:main}, and that $(\rho(B_m)-\rho(S))m^3\to2\pi^2\gamma(\omega)\sum_k\omega_kk^2$ for the corresponding set $S=\{0\}\cup[2,m]$. Computation agrees: for the weights $(\omega_2,\omega_3)=(1,0.4)$, $(\omega_2,\omega_3,\omega_4)=(0.5,1,0.5)$, $(\omega_3,\omega_4,\omega_5)=(0.3,1,0.3)$ and $(\omega_2,\omega_5)=(1,0.7)$, among the sets with holes at depth at most $6$ at one end the satellite minimizes $(\rho(B_m)-\rho(S))m^3$ for $m=40$, $80$, $160$, and the extrapolated values agree within $1\%$ with $2\pi^2\gamma(\omega)\sum_k\omega_kk^2=69.9$, $152.3$, $169.3$ and $104.5$.
\end{rem}

The following constants for the block are used in Sections~\ref{sec:holes} and~\ref{sec:exact}.

\begin{lem}\label{lem:blockconst}
\begin{enumerate}
\item For $n\ge b$, $W-\rho(B_n)\le\bigl(\frac{\pi}{n+1}\bigr)^2\bigl(\sum_{k\in D}\omega_kk^2+\frac{1}{3(n+1)}\sum_{k\in D}\omega_k(k+1)^3\bigr)$; in particular $W-\rho(B_n)\le C_1n^{-2}$ with $C_1=\pi^2\sum_k\omega_k\bigl(k^2+\frac16(k+1)^3\bigr)$.
\item There are $n_B$ and $c_B>0$ with $\rho(B_n)-\rho(B_{n-1})\ge c_Bn^{-3}$ for every $n\ge n_B$.
\item For $\omega\equiv1$ and $D=\{2,3\}$, $\{2,3,4\}$, $\{3,4,5\}$, the bound in~(2) holds for every $n\ge2b$ with $c_B$ as in Table~\ref{tab:const}.
\end{enumerate}
\end{lem}

\begin{proof}
(1) Put $\theta=\pi/(n+1)$ and $x_i=\sin(\theta(i+1))$ for $0\le i<n$, so $\|x\|^2=(n+1)/2$. For $k\in D$, the identity $\sum_{j=1}^{n-k}\sin j\theta\,\sin(j+k)\theta=\frac{n-k}2\cos k\theta+\frac{\sin(k+1)\theta}{2\sin\theta}$ follows from the product formula and $\sin(n-k)\theta=\sin(k+1)\theta$. Hence
\[
W-\frac{x^\top A_{B_n}x}{\|x\|^2}=\frac{1}{n+1}\sum_{k\in D}\omega_k\Bigl(2(n-k)(1-\cos k\theta)+2\Bigl(k+1-\frac{\sin(k+1)\theta}{\sin\theta}\Bigr)\Bigr).
\]
Using $1-\cos k\theta\le k^2\theta^2/2$, $\sin(k+1)\theta\ge(k+1)\theta-(k+1)^3\theta^3/6$ and $\sin\theta\le\theta$ gives the bound.

(2) Apply the proof of Proposition~\ref{prop:ends} with $\mathcal P_J$ replaced by $\{\emptyset\}$, so that the constants $\vartheta_*$, $g_*$, $C$, $c_0$, $\ell$ refer to $P=\emptyset$ only. Put $\Gamma_n(\vartheta)=\vartheta(n-1)+2\Phi_\emptyset(\vartheta)$, the quantization function of $B_n=[0,n-1]$. For large $n$ it is strictly increasing on $[0,\vartheta_*]$ with a unique root $\vartheta_n$ of $\Gamma_n=\pi$, and Lemma~\ref{lem:approx} with the positivity argument there gives $|\rho(B_n)-f(\vartheta_n)|\le C'r_2^{\,n-1}$ for a constant $C'$. Since $\Gamma_{n-1}=\Gamma_n-\vartheta$, we get $\Gamma_n(\vartheta_{n-1})=\pi+\vartheta_{n-1}>\pi$, so $\vartheta_{n-1}>\vartheta_n$, and
\[
\vartheta_n<\vartheta_{n-1}=\Gamma_n(\vartheta_{n-1})-\Gamma_n(\vartheta_n)\le(n-1+2\ell)(\vartheta_{n-1}-\vartheta_n).
\]
As in that proof, $f(\vartheta_n)-f(\vartheta_{n-1})\ge\frac4\pi\sum_k\omega_kk^2\,\vartheta_n(\vartheta_{n-1}-\vartheta_n)\ge\frac4\pi\sum_k\omega_kk^2\,\vartheta_n^2/(n-1+2\ell)$, and $\vartheta_n\ge\pi/(n-1+2g_*)$. Hence $\rho(B_n)-\rho(B_{n-1})\ge c\,n^{-3}-2C'r_2^{\,n-2}\ge\frac c2n^{-3}$ for large $n$.

(3) For $2b\le n\le301$ the spectral radii were enclosed by the Collatz--Wielandt inequalities $\min_i(Av)_i/v_i\le\rho\le\max_i(Av)_i/v_i$, evaluated in exact rational arithmetic for a rational approximation $v>0$ of the Perron vector. For $n\ge301$ we argue as in~(2), now with the certified enclosures. With $\Gamma_n(\vartheta)=\vartheta(n-1)+2\Phi_\emptyset(\vartheta)$ we have $\Gamma_{n-1}=\Gamma_n-\vartheta$. Let $\vartheta_n$ be the smallest root of $\Gamma_n=\pi$; then $\pi/(n-1+2g_{\mathrm{hi}})\le\vartheta_n\le\pi/(n-1)<\theta_{\max}$, where $[g_{\mathrm{lo}},g_{\mathrm{hi}}]$ encloses $g_\emptyset$ on $[0,\theta_{\max}]$. Since $\arctan$ is $1$-Lipschitz, $\Gamma_n(\vartheta)-\Gamma_n(\vartheta_n)\le(\vartheta-\vartheta_n)(n-1+2g_{\mathrm{hi}})+2\vartheta_n(g_{\mathrm{hi}}-g_{\mathrm{lo}})$ for $\vartheta\ge\vartheta_n$, so with $\tau=\vartheta_n/(2(n-1+2g_{\mathrm{hi}}))$ and $2(g_{\mathrm{hi}}-g_{\mathrm{lo}})<\frac12$ every root of $\Gamma_{n-1}=\pi$ exceeds $\vartheta_n+\tau$. As $|f'(\vartheta)|\ge\frac4\pi\sum k^2\,\vartheta$ for $b\vartheta\le\pi/2$,
\[
\rho(B_n)-\rho(B_{n-1})\ge f(\vartheta_n)-f(\vartheta_n+\tau)-2\delta\ge\frac{2\pi\sum_{k\in D}k^2}{(n-1+2g_{\mathrm{hi}})^3}-2\delta,
\]
where $\delta=2|D|\,Cr^{\,n-1+b}/y_{\min}$ bounds the error of Lemma~\ref{lem:approx} for both $B_n$ and $B_{n-1}$. The positivity of $y$ and the enclosures of $g_\emptyset$, $r$, $C$ and $\psi_\emptyset$ near the end were certified by ball arithmetic on a cover of $[0,\theta_{\max}]$ by $200$ balls: there is an integer $J_c$ (Table~\ref{tab:const}) such that for $j\ge J_c$ positivity follows from $A_\emptyset\sin(\vartheta j+\Phi_\emptyset)\ge\frac2\pi\min(j,e-j)$ and $C r^{\,j+b}<\frac1\pi j$, and for $j<J_c$ from the enclosures. Since $n^3r^{\,n}$ decreases for $n\ge301$, the resulting bound at $n=301$ holds for all larger $n$.
\end{proof}

\begin{table}[ht]
\centering
\begin{tabular}{ccccc}
\toprule
$D$ & $c_B$ & $g_\emptyset$ on $[0,\theta_{\max}]$ & $r\le$ & $J_c$\\
\midrule
$\{2,3\}$ & $79.78$ & $[1.68343,1.68370]$ & $0.581$ & $2$\\
$\{2,3,4\}$ & $177.03$ & $[1.95270,1.95313]$ & $0.604$ & $2$\\
$\{3,4,5\}$ & $302.71$ & $[2.37319,2.37411]$ & $0.709$ & $3$\\
\bottomrule
\end{tabular}
\caption{Certified constants for the block; $c_B$ is the minimum of the values for $2b\le n\le301$ and for $n\ge301$.}\label{tab:const}
\end{table}

\subsection{Holes away from the ends}\label{sec:holes}

Throughout this subsection let $S$ be an extremal set for $\rho_m(\omega)$, $\lambda=\rho(S)$, $x$ its unit Perron vector extended by $0$ to $\Z$, $\varepsilon=W-\lambda$, $\mu=\min_Sx$, $L=\min S$, $R=\max S$, and $\Pi(z)=(Ax)(z)=\sum_{y\in S}\omega_{y-z}x_y$ for $z\in\Z$, the \emph{potential} of $x$ at $z$. We first record the Perron transfer, which controls where holes can occur.

\begin{prop}\label{prop:hole}
Let $S$ be an extremal set for $\rho_m(\omega)$, and let $\lambda$, $x$ and $\Pi$ be as above. Then
\[
\Pi(w)\ \le\ (\lambda+\omega_{\max})\min_{v\in S}x_v\qquad\text{for every }w\notin S.
\]
\end{prop}

\begin{proof}
Let $v\in S$ attain $\min_Sx$, let $w\notin S$, and put $S'=(S\setminus\{v\})\cup\{w\}$. Define $x'$ on $S'$ by $x'_w=x_v$ and $x'_y=x_y$ for $y\ne w$. Then $\|x'\|=\|x\|$, and since $\sum_{y\in S}\omega_{y-v}x_y=\lambda x_v$,
\[
x'^\top A_{S'}x'=x^\top A_Sx-2x_v\lambda x_v+2x_v\sum_{y\in S\setminus\{v\}}\omega_{y-w}x_y.
\]
If $\Pi(w)>(\lambda+\omega_{\max})x_v$, then the last sum is at least $\Pi(w)-\omega_{\max}x_v>\lambda x_v$, so the Rayleigh quotient of $x'$ exceeds $\lambda$ and $\rho(S')>\rho(S)$, contradicting extremality.
\end{proof}

On $S$ itself $\Pi=\lambda x\ge\lambda\min_Sx$, so Proposition~\ref{prop:hole} says that every hole sits where the potential $\Pi$ is within a factor $1+\omega_{\max}/\lambda$ of its minimum over $S$. We use the following facts about an extremal set.
\begin{enumerate}
\item[(F1)] $\Pi(z)\le M_0:=(\lambda+\omega_{\max})\mu$ for every $z\notin S$ (Proposition~\ref{prop:hole}).
\item[(F2)] $S$ induces a connected subgraph: otherwise $\rho(S)=\rho(C)$ for a component $C$ with $|C|<m$, and adding $m-|C|$ points one at a time, each adjacent to the current set, gives a connected set of size $m$ with larger spectral radius by the Perron--Frobenius theorem. Since every edge has length at most $b$, $S$ omits no $b$ consecutive integers between $L$ and $R$.
\item[(F3)] $\sum_{z\notin S}\Pi(z)=\varepsilon\sum_Sx$, because $\sum_{z\in\Z}\Pi(z)=W\sum_Sx$ and $\Pi=\lambda x$ on $S$. The points $L-k$ and $R+k$, $k\in D$, lie outside $S$ and have $\Pi(L-k)\ge\omega_kx_L$ and $\Pi(R+k)\ge\omega_kx_R$, so $W\mu\le\varepsilon\sum_Sx\le\varepsilon\sqrt m$.
\item[(F4)] If $m>b$, then $\omega_{\max}\le\rho(B_m)\le\lambda\le W$ and $\varepsilon\le C_1m^{-2}$ by Lemma~\ref{lem:blockconst}(1).
\end{enumerate}
Let $s_D$ be the largest graph distance in $G(\Z,D)$, the graph of the edges of positive weight, from $0$ to an integer $r$ with $0<|r|<b$, minus one; it is finite because $G(\Z,D)$ is connected, and $s_D=1$ for the three sets of Theorem~\ref{thm:all}(2).

\begin{lem}\label{lem:near}
Let $m>b$, put $\Upsilon=\omega_{\min}^{-1}(\lambda/\omega_{\min})^{s_D}$, and let $w\notin S$ and $u\in S$ with $|u-w|\le b$. Then $x_u\le\Upsilon\,\Pi(z)$ for some $z\notin S$ with $|z-w|\le(s_D+2)b$. In particular $x_u\le\Upsilon M_0$, and $\Pi(u)\le\lambda\Upsilon M_0\le(W/\omega_{\min})^{s_D+1}M_0$.
\end{lem}

\begin{proof}
The graph distance from $u$ to $w$ is at most $s_D+1$: by definition of $s_D$ if $|u-w|<b$, and it equals $1$ if $|u-w|=b$. Let $u=y_0,y_1,\dots,y_p=w$ be a shortest path, $p\le s_D+1$. If $y_1,\dots,y_{p-1}\in S$, the eigenvalue equation at $y_{i+1}$ gives $\omega_{\min}x_{y_i}\le\lambda x_{y_{i+1}}$, and $\omega_{\min}x_{y_{p-1}}\le\Pi(w)$, so $x_u\le\omega_{\min}^{-1}(\lambda/\omega_{\min})^{p-1}\Pi(w)$. Otherwise let $y_q$ be the first vertex outside $S$; then likewise $x_u\le\omega_{\min}^{-1}(\lambda/\omega_{\min})^{q-1}\Pi(y_q)$. Since $\lambda\ge\omega_{\max}\ge\omega_{\min}$ by~(F4), $p-1\le s_D$ and each $y_i$ is within $pb\le(s_D+1)b$ of $u$, the first claim follows, and the second by~(F1). Finally $\Pi(u)=\lambda x_u$.
\end{proof}

The flux argument is the following. The function $\Pi$ satisfies $A\Pi=\lambda\Pi+q$, where $q$ is supported near the holes and the ends. Summing this identity from the left end up to a hole $w$ shows that the accumulated flux grows linearly in $w-L$, with a slope of order $\mu$ coming from the left end; but the same flux is a weighted sum of $\Pi$ on a window around $w$, where $\Pi$ is of order $\mu$ by Lemma~\ref{lem:near}. Hence $w-L$ is bounded. Formally, put $q(z)=\sum_{y\notin S}\omega_{y-z}\Pi(y)$. Since $A\Pi=A(\Pi\mathbf 1_S)+A(\Pi\mathbf 1_{\Z\setminus S})=\lambda\Pi+q$, the function
\[
\Xi(t)=\sum_{z<t}\bigl(q(z)-\varepsilon\Pi(z)\bigr)
\]
satisfies $\Xi(t+1)-\Xi(t)=(A\Pi)(t)-W\Pi(t)$. Put $Y(t)=\sum_{s'<t}\Xi(s')$ and $K(n)=\sum_{k\in D}\omega_k(k-|n|)_+$. Then
\begin{equation}\label{eq:Y}
Y(t)=\sum_{z\in\Z}K(t-1-z)\,\Pi(z).
\end{equation}
Indeed, both sides vanish for $t<L-3b$, and both have second difference $(A\Pi)(t-1)-W\Pi(t-1)$ at $t$, because the second difference of $(k-|n|)_+$ is $\mathbf 1[n=k]+\mathbf 1[n=-k]-2\cdot\mathbf 1[n=0]$. Consequently $Y(t+1)\ge K_0\Pi(t)$ with $K_0=K(0)=\sum_k\omega_kk$, and $\sum_nK(n)=\sum_k\omega_kk^2$. Also $\Xi(t)=0$ for $t\le L-2b$, since $\Pi$ and $q$ vanish below $L-b$ and $L-2b$ respectively.

\begin{prop}\label{prop:bulk}
Put
\[
j_1=4b+\frac{4\,(W/\omega_{\min})^{s_D+1}(W+\omega_{\max})\sum_k\omega_kk^2}{W^2},\qquad c''=\Bigl(\frac{WK_0}{2C_1(W+\omega_{\max})(2b+j_1)}\Bigr)^{1/2}.
\]
Then $c''<1$, since $C_1\ge\pi^2\sum_k\omega_kk^2>K_0$ by its definition in Lemma~\ref{lem:blockconst}(1). In the unweighted case $j_1=4b+(2|D|)^{s_D+1}(2|D|+1)\sum_kk^2/|D|^2$.
\begin{enumerate}
\item No hole $w$ satisfies $j_1<w-L$ and $w-L+2b+1\le c''m$; symmetrically from the right end.
\item There is $M_1(\omega)$ such that for $m\ge M_1$ every hole of an extremal set of size $m$ has depth at most $j_1$.
\end{enumerate}
\end{prop}

\begin{proof}
(1) Suppose such a hole exists. Then $R-L>j_1$, so $m>j_1/b$ by~(F2); since $W\ge\omega_{\min}$, $W+\omega_{\max}\ge W$ and $\sum_k\omega_kk^2\ge\omega_{\min}b^2$, we have $j_1\ge4b+4b^2$, so $m>b$ and (F4) applies. Let $w$ be the smallest hole with $w-L>j_1$; then $J':=w-L+2b+1\le c''m$. Every $y\notin S$ with $y<w+b$ and $\Pi(y)>0$ lies in $[L-b,L-1]$, in $(L,L+j_1]$, or in $[w,w+b)$, so there are at most $2b+j_1$ of them, each with $\Pi(y)\le M_0$ by~(F1). For $s'\le w$, every $y$ that contributes to $\sum_{z<s'}q(z)$ satisfies $y<s'+b\le w+b$, and each contributes at most $W\Pi(y)$; hence
\[
\Xi(s')\le\sum_{z<s'}q(z)\le WM_0(2b+j_1)\qquad(s'\le w),
\]
and $Y(t)\le WM_0(2b+j_1)(t-L+2b)$ for $t\le w+1$. With $\Pi(z)\le Y(z+1)/K_0$ we get, for $t\le w$,
\[
\varepsilon\sum_{z<t}\Pi(z)\le\frac{\varepsilon\,WM_0(2b+j_1)}{K_0}\sum_{z=L-b}^{t-1}(z+1-L+2b)\le\frac{\varepsilon WM_0(2b+j_1)J'^2}{2K_0}\le\frac{W^2\mu}4,
\]
using $\varepsilon\le C_1m^{-2}$, $M_0\le(W+\omega_{\max})\mu$ and $J'\le c''m$. On the other side, for $s'\ge L+b$ the points $L-k$ ($k\in D$) lie outside $S$, have $\Pi(L-k)\ge\omega_kx_L\ge\omega_k\mu$, and all the points $z$ with $\omega_{z-(L-k)}>0$ lie below $s'$, so $\sum_{z<s'}q(z)\ge\sum_{k\in D}W\omega_k\mu=W^2\mu/2$. Hence $\Xi(s')\ge W^2\mu/4$ for $L+b\le s'\le w$, and $\Xi(s')\ge-W^2\mu/4$ for $L-2b<s'<L+b$, while $\Xi(s')=0$ for $s'\le L-2b$. Summing, $Y(w)\ge\frac{W^2\mu}4(w-L-4b)$. On the other hand, by~\eqref{eq:Y} $Y(w)$ involves only $z$ with $|z-(w-1)|<b$, hence $|z-w|\le b$, and there $\Pi(z)\le(W/\omega_{\min})^{s_D+1}M_0\le(W/\omega_{\min})^{s_D+1}(W+\omega_{\max})\mu$ by Lemma~\ref{lem:near} and~(F1). Thus $Y(w)\le\sum_k\omega_kk^2\,(W/\omega_{\min})^{s_D+1}(W+\omega_{\max})\mu$, and comparing gives $w-L\le j_1$, a contradiction.

(2) Call a hole \emph{middle} if its depth exceeds $j_1$. By~(1), applied from both ends, a middle hole $w$ satisfies $w-L>c''m-2b-1$ and $R-w>c''m-2b-1$. Suppose middle holes exist, let $w_1<\dots<w_r$ be all of them, and cut $S$ into the pieces $P_0=S\cap[L,w_1)$, $P_i=S\cap(w_i,w_{i+1})$ and $P_r=S\cap(w_r,R]$. Let $m$ be so large that $c''m-2b-1>j_1$. Then every hole inside a piece is within $j_1$ of $L$ or of $R$, so inner pieces are intervals and $P_0$, $P_r$ have at most $j_1$ holes in their hulls. Each piece misses all points of $S$ on one side of some middle hole, and there are at least $\lfloor(c''m-2b-1)/b\rfloor$ of them by~(F2). Hence each piece lies in an interval of at most
\[
n^*=m-\Bigl\lfloor\frac{c''m-2b-1}{b}\Bigr\rfloor+\lfloor j_1\rfloor
\]
points, and $\rho(P_i)\le\rho(B_{n^*})$.

Write $x=\sum_ix_{P_i}$. Every pair $u<v$ of points in different pieces with $\omega_{v-u}>0$ satisfies $u<w<v$ for a middle hole $w$ and $v-u\le b$, so
\[
\lambda=x^\top A_Sx\le\max_i\rho(P_i)+\omega_{\max}\sum_{w}\Bigl(\sum_{u\in S,\ |u-w|<b}x_u\Bigr)^2,
\]
the sum over middle holes $w$. By Lemma~\ref{lem:near}, each $x_u$ in the inner sum is at most $\Upsilon\Pi(z)$ for some $z\notin S$ with $|z-w|\le(s_D+2)b$, and also at most $\Upsilon M_0$. A given $z$ is within $(s_D+2)b$ of at most $2(s_D+2)b+1$ middle holes, so with~(F3),
\[
\sum_w\Bigl(\sum_{|u-w|<b}x_u\Bigr)^2\le2b\Upsilon M_0\cdot(2b)\Upsilon(2(s_D+2)b+1)\sum_{z\notin S}\Pi(z)\le C_X\,\mu\,\varepsilon\sqrt m\le\frac{C_X'}{m^3},
\]
where $C_X$, $C_X'$ depend only on $\omega$, by~(F3), (F4) and $\Upsilon\le\omega_{\min}^{-1}(W/\omega_{\min})^{s_D}$. Therefore
\[
\rho(B_m)\le\lambda\le\rho(B_{n^*})+\frac{C_X'}{m^3}\le\rho(B_m)-(m-n^*)\frac{c_B}{m^3}+\frac{C_X'}{m^3}
\]
by Lemma~\ref{lem:blockconst}(2), provided $n^*\ge n_B$. Since $m-n^*\ge c''m/b-j_1-4$, this is impossible for large $m$.
\end{proof}

\begin{proof}[Proof of Theorem~\ref{thm:all}(1)]
Let $m\ge M_1(\omega)$ and let $S$ be extremal of size $m$. By Proposition~\ref{prop:bulk}, every hole of $S$ has depth at most $j_1$, and $S$ induces a connected subgraph by~(F2). By Proposition~\ref{prop:ends} with $J=\lceil j_1\rceil$, if also $m\ge M(\omega,\lceil j_1\rceil)$, then $\rho(S)\le\rho(B_m)$ with equality only for translates of $B_m$. Since $S$ is extremal, $S$ is a translate of $B_m$.
\end{proof}

\subsection{Exact thresholds}\label{sec:exact}

For the rest of the paper $\omega\equiv1$ and $D$ is one of $\{2,3\}$, $\{2,3,4\}$, $\{3,4,5\}$, so that $W=2|D|$ and $\omega_{\max}=1$ in the results of Section~\ref{sec:def} and Sections~\ref{sec:half}--\ref{sec:holes}; we write $D=[a,b]$. The proof of Theorem~\ref{thm:all}(2) is an induction on $m$ with three ingredients. Under the induction hypothesis the minimum entry $\mu$ of the Perron vector is bounded below, which bounds the edge boundary $|\partial S|=\#\{(z,y):z\notin S,\ y\in S,\ |z-y|\in D\}$ of an extremal set and hence its number of holes (Lemma~\ref{lem:budget}). A hole at depth $d\ge2$ is excluded by playing the two sides of Proposition~\ref{prop:hole} against each other: the Perron vector at the nearer end is bounded above through a resolvent of the stretch between the end and the hole, and the potential at the hole is bounded below through the Perron vector just beyond it (Lemma~\ref{lem:exclude}). Finally, the sets whose only holes have depth one are compared with the block through the approximate eigenvectors of Section~\ref{sec:ends} (Proposition~\ref{prop:satellite}). Each step reduces to inequalities between quantities that are computed with certified enclosures, as described at the end of the subsection.

For a finite $Y\subseteq\Z$, a point $z\notin Y$ and $\nu>\rho(A_Y)$, write
\[
G^Y_\nu(y,z)=\bigl((\nu I-A_Y)^{-1}A_{Y,\{z\}}\bigr)(y)\qquad(y\in Y);
\]
for $\nu=2|D|$ this is the probability that the walk from $y$ leaves $Y$ at $z$. The resolvent $(\nu I-A_Y)^{-1}=\sum_{n\ge0}\nu^{-n-1}A_Y^n$ is entrywise increasing in $Y$ and decreasing in $\nu$. Throughout, $S$ is an extremal set of size $m$, $L=\min S$, and $\lambda_{\mathrm{lo}}$ is any number with $\lambda_{\mathrm{lo}}\le\rho(B_m)$: either a certified Collatz--Wielandt lower bound or $\lambda_{\min}(m)=2|D|-(\frac{\pi}{m+1})^2(\sum k^2+\frac1{3(m+1)}\sum(k+1)^3)$ from Lemma~\ref{lem:blockconst}(1).

\begin{lem}\label{lem:budget}
Suppose $\rho_{m-1}(D)=\rho(B_{m-1})$ and put $g_m=\rho(B_m)-\rho(B_{m-1})$. Then
\[
\mu^2\ge\frac{g_m}{2|D|+g_m}\qquad\text{and}\qquad|\partial S|\le K_{\mathrm{inc}}(m):=\bigl(2|D|-\rho(B_m)\bigr)\sqrt m\,\sqrt{\frac{2|D|+g_m}{g_m}}.
\]
\end{lem}

\begin{proof}
Let $\mu=x_v$. The vector $x-x_ve_v$ restricted to $S\setminus\{v\}$ has norm $\sqrt{1-\mu^2}$ and quadratic form $\lambda-2\lambda\mu^2$, so $\rho(B_{m-1})\ge\rho(S\setminus\{v\})\ge\lambda-\lambda\mu^2/(1-\mu^2)$. With $\rho(B_m)\le\lambda\le2|D|$ this gives $\mu^2/(1-\mu^2)\ge g_m/\lambda\ge g_m/(2|D|)$, which is the first inequality. Each $y\in S$ has $2|D|$ neighbors in $\Z$, so $\sum_{z\in\Z}\Pi(z)=2|D|\sum_Sx$, while $\Pi=\lambda x$ on $S$; hence $\sum_{z\notin S}\Pi(z)=(2|D|-\lambda)\sum_Sx\le(2|D|-\rho(B_m))\sqrt m$. The left side is at least $\mu|\partial S|$.
\end{proof}

The quantity $K_{\mathrm{inc}}(m)$ increases to $\pi\sqrt{|D|\sum_{k\in D}k^2}$, which is $16.02$, $29.30$ and $38.48$ for the three sets, and for $m\ge2b$ Lemma~\ref{lem:blockconst} gives $K_{\mathrm{inc}}(m)\le K_u(m):=\pi^2\bigl(\sum k^2+\frac1{3m}\sum(k+1)^3\bigr)\bigl(2|D|/c_B+m^{-3}\bigr)^{1/2}$. Let $h(K)$ be the largest number of holes of a finite set with no $b$ consecutive gaps between its minimum and maximum and with $|\partial S|\le K$. Since $|\partial S|$ is a sum of contributions of pairs at distance at most $b$, $h(K)$ is the length of a longest path in a graph whose vertices record the membership pattern of the last $b$ positions and the boundary so far. A step that adds nothing to the boundary forces $s_i=s_{i-k}$ for every $k\in D$, so a cycle of such steps is periodic with periods $a$ and $a+1$, hence constant, and since $b$ consecutive gaps are excluded it is the all-in pattern, which carries no holes. Hence $h(K)$ is computed by a dynamic program over $(K+1)2^b$ positions. Under the hypothesis of Lemma~\ref{lem:budget}, an extremal set satisfies $\max S-\min S\le m-1+H(m)$ with $H(m)=h(\lfloor K_{\mathrm{inc}}(m)\rfloor)$.

\begin{lem}\label{lem:exclude}
Let $w$ be a hole of $S$.
\begin{enumerate}
\item Put $c_r=\frac12$ for $\{2,3\}$, $(c_1,c_2,c_3)=(\frac12,\frac7{12},\frac12)$ for $\{2,3,4\}$, and $c_r=\frac13$ for $\{3,4,5\}$. Put $d:=w-L$ and suppose $\rho(B_d)<\lambda_{\mathrm{lo}}$. If
\[
(2|D|+1)\,G^{[0,d)}_{\lambda_{\mathrm{lo}}}(0,d+r)<c_r\qquad\text{for }r=1,\dots,b-1,
\]
then $S$ is not extremal.
\item Suppose $j:=w-L\ge1$, let $\Lambda\ge b$ and $M=j+1+\Lambda$, and suppose $(w,L+M+b)\subseteq S$. Let $Y$ be a set with $(S-L)\cap[0,M)\subseteq Y\subseteq[0,M)\setminus\{j\}$ and $\rho(A_Y)<\lambda_{\mathrm{lo}}$, and put
\[
a_r=\sum_{k\in D}G^{(0,\Lambda+1)}_{2|D|}(k,\Lambda+1+r),\qquad W_r=G^{Y}_{\lambda_{\mathrm{lo}}}(0,M+r)\qquad(0\le r<b).
\]
If $a_r>(2|D|+1)W_r$ for every $r$, then $S$ is not extremal.
\end{enumerate}
\end{lem}

\begin{proof}
Translate so that $L=0$. In both parts we show $\max_{z\notin S}\Pi(z)>(2|D|+1)x_0\ge(\lambda+1)\mu$, contradicting Proposition~\ref{prop:hole}.

(1) Put $T=\{r\in(0,b):w+r\in S\}$ and $t_r=x_{w+r}$. On $U=S\cap[0,w)$ the eigenvalue equation reads $(\lambda I-A_U)x_U=A_{U,w+T}\,t$, because every neighbor of $U$ in $S\setminus U$ lies in $w+T$. By the monotonicity of the resolvent and $\lambda\ge\rho(B_m)\ge\lambda_{\mathrm{lo}}$,
\[
x_0\le\sum_{r\in T}t_r\,G^{[0,d)}_{\lambda_{\mathrm{lo}}}(0,d+r).
\]
On the other side, for $y\in S$ the eigenvalue equation and $\lambda\le2|D|$ give
\[
x_y\ge\frac1{2|D|}\sum_{r\in T,\ |w+r-y|\in D}t_r.
\]
Hence $\Pi(z)\ge\sum_{r\in T}c_{z,r}t_r$ for every $z\notin S$, where
\[
c_{z,r}=\mathbf 1\bigl[|w+r-z|\in D\bigr]+\frac{1}{2|D|}\,\#\bigl\{y\in S\setminus(w+T):\ |y-w|\le2b,\ |y-z|\in D,\ |w+r-y|\in D\bigr\}
\]
depends only on the configuration $S\cap[w-2b,w+2b]$. For each of the finitely many configurations, a linear program gives weights $\alpha_z\ge0$ on the points $z\notin S$ of the window with $\sum_z\alpha_z=1$ and $\sum_z\alpha_zc_{z,r}\ge c_r$ for all $r\in T$. Therefore $\max_{z\notin S}\Pi(z)\ge\sum_z\alpha_z\Pi(z)\ge\sum_{r\in T}c_rt_r>(2|D|+1)\,x_0$.

(2) Put $t_r=x_{M+r}$; these points lie in $S$. Let $U=S\cap[0,M)$. Every neighbor of $U$ in $S\setminus U$ lies in $[M,M+b)$, so $(\lambda I-A_U)x_U=A_{U,[M,M+b)}t$, and since $U\subseteq Y$ and $\lambda\ge\lambda_{\mathrm{lo}}$, monotonicity of the resolvent gives $x_0\le\sum_rt_rW_r$. On $U'=(j,M)\subseteq S$ the eigenvalue equation reads $(\lambda I-A_{U'})x_{U'}=A_{U',[M,M+b)}t+q$, where $q\ge0$ collects the values of $x$ at the points of $S$ below $j$. Since $\rho(A_{U'})<\lambda\le2|D|$, this gives $x_{U'}\ge(2|D|I-A_{U'})^{-1}A_{U',[M,M+b)}t$, and after translating $j$ to $0$,
\[
\Pi(j)\ge\sum_{k\in D}x_{j+k}\ge\sum_rt_ra_r>(2|D|+1)\sum_rt_rW_r\ge(2|D|+1)\,x_0.\qedhere
\]
\end{proof}

For a given $m$ the resolvents in Lemma~\ref{lem:exclude} are computed in ball arithmetic, and the positivity of the computed column certifies $\lambda_{\mathrm{lo}}>\rho(B_d)$, respectively $\lambda_{\mathrm{lo}}>\rho(A_Y)$. In part~(2) we take $Y=[0,M)\setminus\{j\}$. At the depth $j=2$ for $\{3,4,5\}$ this does not suffice, and we use the same argument with more information. For each of the two possible patterns $P\in\{\{0\},\{0,1\}\}$ of $S$ below the hole, $U=S\cap[0,M)=P\cup(2,M)$ is known, so the eigenvalue equation gives $x_U\ge(2|D|I-A_U)^{-1}A_{U,[M,M+b)}t$ on all of $U$ and $x_y\le\sum_rt_rG^U_{\lambda_{\mathrm{lo}}}(y,M+r)$ for every $y\in U$. We then compare $\Pi(z)$ with $(2|D|+1)x_y$ for a point $z\in[-b,2]\setminus S$ and a point $y\in U$ chosen to make the inequality strongest, instead of $\Pi(j)$ with $(2|D|+1)x_0$; a contradiction again follows from Proposition~\ref{prop:hole}. For large $m$ we compare with the walk instead. If $Y\subseteq[0,M)$, then $\chi(y)=(y+b)(M+b-y)/\sigma^2$ satisfies $(2|D|I-A)\chi=2|D|$, because $\chi(y+k)+\chi(y-k)=2\chi(y)-2k^2/\sigma^2$, and $\chi\ge0$ on $[-b,M+b)$, which contains every exit point, so the exit time satisfies $\E_y[\tau_Y]\le(M+2b)^2/(4\sigma^2)$. Let $\nu=2|D|-\varepsilon$, let $\zeta\colon\Z\setminus Y\to[0,1]$, and put $V=(\nu I-A_Y)^{-1}A_{Y,\Z\setminus Y}\zeta$ and $v=(2|D|I-A_Y)^{-1}A_{Y,\Z\setminus Y}\zeta$. With $G=(2|D|I-A_Y)^{-1}$ we have $G\mathbf 1=\E_\cdot[\tau_Y]/(2|D|)$. If $\varkappa:=\varepsilon\max_{y\in Y}\E_y[\tau_Y]/(2|D|)<1$, then $\|\varepsilon G\|_\infty=\varkappa<1$, so $\nu>2|D|-1/\|G\|_\infty\ge\rho(A_Y)$, and the resolvent identity $V=v+\varepsilon GV$ gives $\|V\|_\infty\le1/(1-\varkappa)$ and
\begin{equation}\label{eq:resolvent}
V(y)\le v(y)+\frac{\varepsilon\,\E_y[\tau_Y]}{2|D|(1-\varkappa)}\qquad(y\in Y).
\end{equation}
With $\varepsilon=2|D|-\lambda_{\mathrm{lo}}$, this bounds $G^{[0,d)}_{\lambda_{\mathrm{lo}}}$ and $W_r$ by the harmonic measures at $\nu=2|D|$ plus a correction of order $\varepsilon\,\E[\tau]$. For the expected exit times we also use $\sigma^2\E_0[\tau_{[0,d)}]\le\min\bigl(b(d+b),\,h_0(0)(d+b)^2/h_0(d)+b^2\bigr)$, from the martingale $X_t^2-\sigma^2t$ and the fact that, $h_0$ being increasing, the walk from $0$ reaches $[n,\infty)$ before $(-\infty,0)$ with probability at most $h_0(0)/h_0(n)$; for $Y=[0,M)\setminus\{j\}$ the extra term $h_0(0)j^2/h_0(j)$ accounts for the exit at $j$. Finally, $v_r(d)=G^{[0,d)}_{2|D|}(0,d+r)\le h_0(0)/h_0(d+r)$.

It remains to compare the block with the sets whose only holes have depth one.

\begin{prop}\label{prop:satellite}
$\rho(\{0\}\cup[2,m])<\rho(B_m)$ and $\rho(\{0\}\cup[2,m-1]\cup\{m+1\})<\rho(B_m)$ for every $m\ge2b+2$.
\end{prop}

\begin{proof}
For $2b+2\le m\le300$ the claim is certified by the Collatz--Wielandt inequalities with a common $\lambda_0$ between the spectral radii. Let $m\ge301$, write $\Phi_{\mathrm{sat}}=\Phi_{\{-2\}}$ and $H=\Phi_\emptyset+\vartheta-\Phi_{\mathrm{sat}}$. The three sets have interval parts of lengths $m-1$, $m-2$, $m-3$, so their quantization functions are $\Gamma_B=\vartheta(m-1)+2\Phi_\emptyset$, $\Gamma_B-H$ and $\Gamma_B-2H$. By ball arithmetic on the cover of $[0,\theta_{\max}]$ used in Lemma~\ref{lem:blockconst}(3), $H/\vartheta\ge H_{\mathrm{lo}}$ with $H_{\mathrm{lo}}$ as in Table~\ref{tab:sat}, and $2(g_{\mathrm{hi}}-g_{\mathrm{lo}})<H_{\mathrm{lo}}/2$. As in that proof, with $\tau=H_{\mathrm{lo}}\vartheta_B/(2(m-1+2g_{\mathrm{hi}}))$, every root of $\Gamma_B-kH=\pi$ ($k=1,2$) in $(0,\theta_{\max}]$ exceeds $\vartheta_B+\tau$, and such a root exists since $\Gamma_B-kH\ge\vartheta(m-3)$ and $\theta_{\max}(m-3)>\pi$. The approximate eigenvectors of Lemma~\ref{lem:approx} are positive by the same argument, now also using the certified enclosure of $\psi_{\mathrm{sat}}(-2)$. Hence
\[
\rho(B_m)-\rho(S)\ge\frac{2\pi H_{\mathrm{lo}}\sum_{k\in D}k^2}{(m-1+2g_{\mathrm{hi}})^3}-\delta_B-\delta_S,
\]
where $\delta_B$, $\delta_S$ are the error bounds of Lemma~\ref{lem:approx} for the two sets. The right side is positive at $m=301$ (Table~\ref{tab:sat}) and stays positive for larger $m$, because the error terms decay like $r^m$.
\end{proof}

\begin{table}[ht]
\centering
\begin{tabular}{cccc}
\toprule
$D$ & $H_{\mathrm{lo}}$ & main term at $m=301$ & $\delta_B+\delta_S\le$\\
\midrule
$\{2,3\}$ & $0.4594$ & $1.34\times10^{-6}$ & $10^{-70}$\\
$\{2,3,4\}$ & $0.3934$ & $2.55\times10^{-6}$ & $10^{-64}$\\
$\{3,4,5\}$ & $0.3231$ & $3.59\times10^{-6}$ & $10^{-43}$\\
\bottomrule
\end{tabular}
\caption{Certified data for Proposition~\ref{prop:satellite}; compare $H_{\mathrm{lo}}$ with $\gamma$ in Remark~\ref{rem:gamma}.}\label{tab:sat}
\end{table}

\begin{proof}[Proof of Theorem~\ref{thm:all}(2)]
We use induction on $m$, in four ranges whose boundaries are listed in Table~\ref{tab:ind}.

\emph{Base} ($m=m_0$). All sets with $\min S=0$, $|S|=m_0$ and no $b$ consecutive gaps are enumerated, and $\rho(S)<\rho(B_{m_0})$ is certified for every connected one except $B_{m_0}$ itself and, for $\{3,4,5\}$, the two sets $[0,10]\setminus\{2\}$ and $[0,10]\setminus\{8\}$; the map fixing $\{0,1,3,\dots,9\}$ and sending $2\mapsto10$ is an isomorphism from $B_{10}$ to $[0,10]\setminus\{2\}$.

\emph{Small $m$} ($m_0<m\le M_e$). By the induction hypothesis and Lemma~\ref{lem:budget}, an extremal set has $|\partial S|\le\lfloor K_{\mathrm{inc}}(m)\rfloor$; the sets with this property and no $b$ consecutive gaps are enumerated, and $\rho(S)<\rho(B_m)$ is certified for each of them other than $B_m$.

\emph{Middle range} ($M_e<m<M_u$). Lemma~\ref{lem:exclude}(1) is checked for every $d$ with $2\le d\le(m+H(m)-1)/2$; let $\mathcal F(m)$ be the set of depths where it fails. Every hole then has depth in $\mathcal F(m)\cup\{1\}$, measured from the nearer end. Let $w$ be the hole nearest to the middle on the left side, suppose that its depth $j$ is at least $2$, and put $Z=\max\mathcal F(m)$. Every hole beyond $w$ lies within $Z$ of $\max S$, so for every $\Lambda$ with $j+\Lambda+b+Z+1\le m-1$ the hypothesis $(w,L+M+b)\subseteq S$ of Lemma~\ref{lem:exclude}(2) holds, and also $\max S-w>j$. Lemma~\ref{lem:exclude}(2) is checked for every $j\in\mathcal F(m)$ with $j\ge2$ and the first $\Lambda\in\{5,10,20,40\}$ satisfying this inequality for which it succeeds; for $\{3,4,5\}$ and $j=2$ the variant described after Lemma~\ref{lem:exclude} is used when this fails. This excludes every such $j$, except that for $\{3,4,5\}$ and $68\le m\le110$ the depth $j=2$ survives; the $15$ sets with holes at depths in $\{1,2\}$ at each end, other than the block, are then compared with $B_m$ directly. By symmetry the same holds at the right end. Hence every hole has depth $1$, so $S$ is $B_m$, $\{0\}\cup[2,m]$ or $\{0\}\cup[2,m-1]\cup\{m+1\}$ up to translation and reflection, and Proposition~\ref{prop:satellite} excludes the last two.

\emph{Large $m$} ($m\ge M_u$). Here $H(m)\le H_u:=h(\lfloor K_u(M_u)\rfloor)$, and Lemma~\ref{lem:exclude}(1) is applied with $\lambda_{\mathrm{lo}}=\lambda_{\min}(m)$ and the bound~\eqref{eq:resolvent}, whose correction term decreases in $m$. For $d\le(M_u+H_u-1)/2$ it therefore suffices to check $m=M_u$. For larger $d$ we have $m\ge2d+1-H_u$, and the check is made at $m=\max(M_u,2d+1-H_u)$ for $d\le600$; beyond that the bound decreases in $d$, since $\varepsilon\E_0[\tau]=O(1/d)$ and $\varkappa$ is proportional to $((d+2b)/(2d+2-H_u))^2$, so the check at $d=601$ covers all larger $d$. The depths at which Lemma~\ref{lem:exclude}(1) fails are listed in Table~\ref{tab:ind}; Lemma~\ref{lem:exclude}(2) with the bound~\eqref{eq:resolvent} excludes them for all $m\ge M_u$, except the depth $j=2$ for $\{3,4,5\}$. Its hypothesis $(w,L+M+b)\subseteq S$ holds as in the middle range: the depths outside this list are excluded by Lemma~\ref{lem:exclude}(1) for every $m\ge M_u$, so with $Z$ the largest depth in the list, every hole beyond $w$ lies within $Z$ of $\max S$, and $j+\Lambda+b+Z+1\le M_u-1\le m-1$. There the variant described after Lemma~\ref{lem:exclude} is checked with $\lambda_{\mathrm{lo}}=\lambda_{\min}(M_u)$; since $\lambda\ge\lambda_{\min}(m)\ge\lambda_{\min}(M_u)$ and the resolvent is decreasing in $\lambda_{\mathrm{lo}}$, this covers all $m\ge M_u$. Proposition~\ref{prop:satellite} handles depth $1$.
\end{proof}

\begin{table}[ht]
\centering
\begin{tabular}{cccccccc}
\toprule
$D$ & $m_0$ & $M_e$ & sets & per-$m$ range & $M_u$ & $\Lambda$ ($m\ge M_u$) & depths left\\
\midrule
$\{2,3\}$ & $7$ & $41$ & $6.4\times10^3$ & $27$--$320$ & $300$ & $5$ & $2,4,5$\\
$\{2,3,4\}$ & $8$ & $71$ & $1.6\times10^5$ & $41$--$470$ & $400$ & $10$ & $2,3,5,6,7$\\
$\{3,4,5\}$ & $10$ & $67$ & $5.9\times10^6$ & $68$--$520$ & $500$ & $10$, $20$ & $2,\dots,13$\\
\bottomrule
\end{tabular}
\caption{Ranges in the proof of Theorem~\ref{thm:all}(2). Sets are enumerated for $m_0\le m\le M_e$, including the base; the per-$m$ computations of the middle range succeed on the range shown, which contains $(M_e,M_u)$; the last column lists the depths at which Lemma~\ref{lem:exclude}(1) fails for $m\ge M_u$, which are then excluded by Lemma~\ref{lem:exclude}(2). For $\{3,4,5\}$, $\Lambda=20$ is used at the exceptional depth $j=2$.}\label{tab:ind}
\end{table}

All computations were done in Python, with ball arithmetic from Arb~\cite{Arb} through python-flint. Spectral radii are enclosed by the Collatz--Wielandt inequalities of Lemma~\ref{lem:approx}, evaluated for an integer rounding $v$ of the Perron vector and checked in exact integer arithmetic; in the enumerations for $\{3,4,5\}$ a floating-point vector is used instead, with a rigorous bound on the rounding error of each residual. Harmonic measures and resolvents are solutions of linear systems with the roots of $\sum_{k\in D}(z^k+z^{-k})=\nu$ as data, and these systems, the half-line functions of Section~\ref{sec:ends} on a cover of $[0,\theta_{\max}]$ by balls, and the resolvents for given $m$ were solved in ball arithmetic. The weights of the linear programs in Lemma~\ref{lem:exclude}(1) were rationalized and the inequalities checked in exact arithmetic. The enumeration of candidate sets was checked against brute force for $m\le9$. The code and its output logs are available as ancillary files with the arXiv version of this paper; a README there maps each lemma and table to the script and log that certify it.

\section{Concluding Remarks}\label{sec:remarks}

Theorem~\ref{thm:all}(1) gives no value of $m_0(\omega)$. The proof is effective in the sense that $m_0(\omega)$ is determined by the constants of Proposition~\ref{prop:bulk}, which are explicit, and by the frequency range on which the phase deficit of each end pattern of depth at most $j_1$ exceeds $\gamma(\omega)\vartheta/2$, which can be computed from $\omega$; we have not extracted a bound. The quantity $j_1$ grows like $(W/\omega_{\min})^{s_D}$, which is exponential in $b$ when $s_D$ is of order $b$, as for $D=\{b-1,b\}$, and the number of patterns grows exponentially in $j_1$, so a bound obtained this way would be far larger than the true threshold.

For the application to reserve design, two features of real connectivity matrices are harmless. A drift along the coast of the form $C_{ij}=\omega_{j-i}e^{\alpha(j-i)}$ is removed from every principal submatrix by the diagonal similarity $C\mapsto\Delta C\Delta^{-1}$ with $\Delta=\mathrm{diag}(e^{\alpha i})$, and a constant self-recruitment rate $C_{ii}$ only shifts the spectral radius. Hence Theorem~\ref{thm:all}(1) covers every distance-dependent dispersal kernel of finite range whose support has greatest common divisor $1$.

\section*{Acknowledgments}

The author thanks Sangil Kim for discussions on larval connectivity matrices. The author used Claude (Anthropic) for language refinement, technical editing, computational verification of examples, and symbolic computations. The author is responsible for all mathematical content of this paper.

\paragraph{Data availability.} The code and output logs that certify the computations in Section~\ref{sec:exact} are available as ancillary files with the arXiv version of this paper.

\end{document}